\documentclass[twocolumn]{autart}    % Comment this line out if you need a4paper

\usepackage{amsfonts}
\usepackage{amsmath}
\usepackage{epstopdf}
\usepackage{graphicx}
\newcommand{\R}{\mathbb{R}}

\newcommand{\bmat}[1]{\begin{bmatrix}#1\end{bmatrix}}

\newcommand{\diag}{\mathop{\mathrm{diag}}}

\newcommand{\mcl}[1]{\mathcal{ #1}}
\newcommand{\mbf}[1]{\mathbf{ #1}}

\newcommand{\fourpi}[4]{\ensuremath{\mcl{P}{\tiny\bmat{#1,& \hspace{-3mm}#2 \\ #3,& \hspace{-1.5mm} \left\{#4\right\} }}}}

\let\bbbl\biggl
\let\bbbbl\Biggl

\let\bbbr\biggr
\let\bbbbr\Biggr

\begin{document}
%
%\title{\LARGE \bf
%Representation of Systems and Networks with Delay:\\ DDEs, DDFs, ODE-PDEs and PIEs
%%Representation of Systems with Delay: Delay-Differential Equations, Differential-Difference Equations, and Partial-Integral Equations
%
%%Convexification of the Synthesis Problem for Infinite-Dimensional Systems with Application to MIMO Multi-Delay Systems and Implementation in SOS
%}

		\begin{frontmatter}
	\title{Representation of Networks and Systems with Delay:\\ DDEs, DDFs, ODE-PDEs and PIEs}

	\author[Peet]{Matthew M. Peet}\ead{mpeet@asu.edu},               % e-mail address
	
%	\address[Paestum]{\textsc{School of Engineering of Matter, Transport and Energy at Arizona State University, Tempe, AZ, 85287-6106 USA }}  % Please supply
	\address[Peet]{\textsc{School for the Engineering of Matter, Transport and Energy, Arizona State University, Tempe, AZ, 85298 USA.}}             % full addresses
\thanks{This work was supported by the National Science Foundation under grant No. 1739990.}

		\begin{keyword}                           % Five to ten keywords,
			Delay, PDEs, Networked Control              % chosen from the IFAC
		\end{keyword}                             % keyword list or with the
		% help of the Automatica
		% keyword wizard

%%%%%%%%%%%%%%%%%%%%%%%%%%%%%%%%%%%%%%%%%%%%%%%%%%%%%%%%%%%%%%%%%%%%%%%%%%%%%%%%
\begin{abstract}
Delay-Differential Equations (DDEs) are the most common representation for systems with delay. However, the DDE representation is limited. In network models with delay, the delayed channels are low-dimensional and accounting for this heterogeneity is not possible in the DDE framework. In addition, DDEs cannot be used to model difference equations. Furthermore, estimation and control of systems in DDE format has proven challenging, despite decades of study. In this paper, we examine alternative representations for systems with delay and provide formulae for conversion between representations. First, we examine the Differential-Difference (DDF) formulation which allows us to represent the low-dimensional nature of delayed information. Next, we examine the coupled ODE-PDE formulation, for which backstepping methods have recently become available. Finally, we consider the algebraic Partial-Integral Equation (PIE) representation, which allows the optimal estimation and control problems to be solved efficiently through the use of recent software packages such as PIETOOLS. In each case, we consider a very general class of delay systems, specifically accounting for all four possible sources of delay - state delay, input delay, output delay, and process delay. We then apply these representations to 3 archetypical network models.
\end{abstract}
		
	\end{frontmatter}
		
\section{Introduction} Delay-Differential Equations (DDEs) are a convenient shorthand notation used to represent what is perhaps the simplest form of spatially-distributed phenomenon - transport. Because of their notational simplicity, it is common to use DDEs to model very complex systems with multiple sources of delay - including almost all models of control over and of ``networks''.

To illustrate the ways in which delays can complicate an otherwise straightforward control problem, consider control of a swarm of $N$ UAVs over a wireless network. In this case, each UAV, $i$, has a state, $x_i(t)\in \R^{n_i}$ which may represent, e.g. displacement (the concatenation of all such states is denoted $x$). Each UAV has local sensors which measure $y_i$ and this information is transmitted to a centralized control authority. There is also a centralized vector of inputs, $u$, a regulated vector of outputs, $z$, and a vector of disturbances, $w$ - including both process and sensor noise. We model this system as follows.
\begin{align}
\dot x_i(t)&=a_ix_i(t)+\sum\nolimits_{j=1}^N a_{ij} x_j(t-\hat \tau_{ij})\notag \\
&\hspace{2cm}+b_{1i}w(t-\bar \tau_i)+b_{2i} u(t-h_i)\notag\\
z(t)&=C_{1} x(t)+D_{12}u(t)\notag \\
y_i(t)&=c_{2i} x_i(t-\tilde \tau_i)+d_{21i}w(t-\tilde \tau_i) \label{eqn:network1}\\[-6mm] \notag
\end{align}

\begin{itemize}
\item $a_i$ is the internal dynamics of the UAV $i$
\item $a_{ij}$ is the effect of UAV $j$ on the state of UAV $i$.
\item $b_{1i}$ is the disturbance to the motion of UAV $i$
\item $b_{2i}$ is the effect of the central command on UAV $i$
\item $c_{2i}$ is the measurement of the state of UAV $i$
\item $d_{21i}$ is the disturbance to the sensor on UAV $i$
\item $C_1$ gives the weight on states of the fleet of UAVs to minimize in the optimal control problem
\item $D_{12}$ gives the weight on actuator commands to minimize in the optimal control problem
\item $\hat \tau_{ij}$ is the time taken for changes in state of UAV $j$ to affect UAV $i$
\item $h_i$ is the time taken for a command from the central authority to reach UAV $i$
\item $\bar \tau_i$ is the time it takes the process disturbance (wind, tracking signal, et c.) to reach UAV $i$
\item $\tilde \tau_i$ is the time taken for measurements collected at UAV $i$ to reach the central authority
\end{itemize}
This relatively simple model shows that delayed channels are often low dimensional ($\R^{n_i}$ vs. $\R^{\sum n_i}$) and specifies four separate yet individually significant sources of delay. Specifically, we have: state delay ($\hat \tau_{ij}$); input delay ($h_i$); process delay ($\bar \tau_i$); and output delay ($\tilde \tau_i$).

This UAV network is modeled as a DDE - a structure formulated in Eqn.~\eqref{eqn:DDE} in Sec.~\ref{sec:DDE}. If we consider control of such a network, however, we find that while there are algorithms for control of DDEs (See~\cite{peet_2020SICON}),  these algorithms are complex and are memory-limited to a relatively small number of UAVs (perhaps 4-5). The premise of this paper, however, is that the limitations of these algorithms are not caused by inefficiency of the algorithms, but rather by the failure to account for the low dimensional nature of the delayed channels. Specifically, we note that in the UAV model, while the concatenated state, $x(t)$, is high-dimensional, the individual delayed channels, $x_i(t)$, are of much lower dimension. If we represent the network as a DDE using the formulation in Subsec.~\ref{subsec:UAV_DDE}, then the low-dimensional nature of the delayed channels is lost. Furthermore, DDEs cannot represent some important system designs - including a model of feedback described in Subsection~\ref{subsec:UAV_SOF}.

For these reasons, in Sec.~\ref{sec:DDF}, we consider the use of Differential Difference Equations (DDFs). The DDF can be used to model both DDEs and neutral-type systems, while also allowing for the assignment of delayed information to heterogeneous low-dimensional channels. Specifically, the infinite-dimensional component of state-space (as defined in~\cite{gu_2010,pepe_2008}) of the UAV network in the DDF framework is $\prod_i L_2[-\tau_i,0]^{n_i}$ as opposed to $\prod_i L_2[-\tau_i,0]^{\sum n_i}$ using a DDE. In addition, DDFs allow us to represent difference equations which arise in some network models - See Subsection~\ref{subsec:UAV_SOF}.

From the DDF model we turn to coupled ODE-PDE models in Sec.~\ref{sec:ODEPDE}. ODE-PDEs can be used to model a variety of systems. However, for the particular class of ODE-PDEs we use in Sec.~\ref{sec:ODEPDE}, the solutions to the ODE-PDE are equivalent to those of the DDF (as defined in Sec.~\ref{sec:DDF}). Backstepping methods have been developed for ODE-PDE models of delay (e.g.~\cite{krstic_2008,zhu_2015}) and the formulae we present for conversion of DDFs to ODE-PDEs may prove useful if the reader is interested in application or further development of these backstepping methods.

% Finally, we consider the
Next, in Sec.~\ref{sec:PIE}, we consider Partial Integral Equations (PIEs)~\cite{appell_book}.
PIEs are a generalization of integro-differential equations of Barbashin type which have been used since the 1950s to model systems in biology, physics, and continuum mechanics (See chapters 19-20 of~\cite{appell_book} for a survey). PIEs and ODE-PDEs define an equivalent set of solutions and in this section, we provide formulae for conversion of DDEs and DDFs to PIEs.
%PIE representations have the form
%\begin{align*}
%\mcl T \dot{\mbf x}(t)+\mcl B_{d1}\dot w(t)+\mcl B_{d2}\dot u(t)&=\mcl A\mbf x(t)+\mcl B_1w(t)+\mcl B_2u(t)\\
%\hspace{-1cm}z(t)=\mcl C_1\mbf x(t)&+\mcl D_{11}w(t)+\mcl D_{12}u(t),\\
%\hspace{-1cm}y(t)=\mcl C_2\mbf x(t)&+\mcl D_{21}w(t)+\mcl D_{22}u(t)
%\end{align*}
%where the $\mcl T, \mcl A, \mcl B_{i}, \mcl C_i, \mcl D_{ij}$ are Partial Integral (PI) operators and have the form
%\[		\left(\fourpi{P}{Q_1}{Q_2}{R_i}\bmat{x\\\mbf \Phi}\right)(s):= {\bmat{
%				Px + \int_{-1}^{0} Q_1(s)\mbf \Phi(s)ds\\
%				Q_2(s)x +\left(\mcl P_{\{R_i\}}\mbf \Phi\right)(s)
%		}}\vspace{-3mm}
%\]
%where
%\begin{align*}
%&\left(\mcl{P}_{\{R_i\}}\mbf \Phi\right)(s):= \\
%&R_0(s) \mbf \Phi(s) +\int_{-1}^s R_1(s,\theta)\mbf \Phi(\theta)d \theta+\int_s^0R_2(s,\theta)\mbf \Phi(\theta)d \theta.
%\end{align*}
PIE models have the advantage that they are defined by Partial Integral (PI) operators. Unlike Dirac and differential operators, PI operators are bounded and form an algebra. Furthermore, PIE models do not require boundary conditions or continuity constraints  - simplifying analysis and optimal control problems. Indeed, it has been recently shown in~\cite{shivakumar_2019CDC,das_2019CDC} that many problems in analysis, optimal estimation and control of ODE-PDE models can be formulated as optimization over the cone of positive PI operators. In Sec.~\ref{sec:temp}, we show that the PIE formulation allows for $H_\infty$-optimal control of a 40 user, 80-state, 40-delay, 40-input, 40-disturbance network model of temperature control.

\begin{figure*}[!t]
\textbf{The Class of Delay-Differential Equations (DDEs):}\vspace{-3mm}
\begin{align}
	&\bmat{\dot{x}(t)\\z(t) \\ y(t)}=\bmat{A_0 & B_1 & B_2\\ C_{10} & D_{11} &D_{12}\\ C_{20} & D_{21} &D_{22}}\bmat{x(t)\\w(t)\\u(t)}+\sum_{i=1}^K \bmat{A_i & B_{1i} & B_{2i}\\C_{1i} & D_{11i} & D_{12i}\\C_{2i} & D_{21i} & D_{22i}} \bmat{x(t-\tau_i)\\w(t-\tau_i)\\u(t-\tau_i)}\notag\\
& \hspace{3cm}+\sum_{i=1}^K \int\limits_{-\tau_i}^0\bmat{A_{di}(s) & B_{1di}(s) &B_{2di}(s)\\C_{1di}(s) & D_{11di}(s) & D_{12di}(s)\\C_{2di}(s) & D_{21di}(s) & D_{22di}(s)} \bmat{x(t+s)\\w(t+s)\\u(t+s)}ds\label{eqn:DDE}
\end{align}
\textbf{The Class of Differential-Difference Equations (DDFs):}\vspace{-3mm}
\begin{align}
	\bmat{\dot{x}(t)\\ z(t)\\y(t)\\r_i(t)}&=\bmat{A_0 & B_1& B_2\\C_1 &D_{11}&D_{12}\\C_2&D_{21}&D_{22}\\C_{ri}&B_{r1i}&B_{r2i}}\bmat{x(t)\\w(t)\\u(t)}+\bmat{B_v\\D_{1v}\\D_{2v}\\D_{rvi}} v(t) \notag\\
v(t)&=\sum_{i=1}^K C_{vi} r_i(t-\tau_i)+\sum_{i=1}^K \int_{-\tau_i}^0C_{vdi}(s) r_i(t+s)ds.\label{eqn:DDF}
\end{align}
\textbf{The Class of Neutral-Type Systems (NDS):}\vspace{-3mm}
\begin{align}
	&\bmat{\dot{x}(t)\\z(t) \\ y(t)}=\bmat{A_0 & B_1 & B_2\\ C_{10} & D_{11} &D_{12}\\ C_{20} & D_{21} &D_{22}}\bmat{x(t)\\w(t)\\u(t)}+\sum_{i=1}^K \bmat{A_i  & B_{1i} & B_{2i}& E_i\\C_{1i}& D_{11i} & D_{12i} & E_{1i}\\C_{2i} & D_{21i} & D_{22i}&E_{2i}} \bmat{x(t-\tau_i)\\w(t-\tau_i)\\u(t-\tau_i)\\ \dot x(t-\tau_i)}\notag\\
&\hspace{2.6cm} +\sum_{i=1}^K \hspace{0mm}\int_{-\tau_i}^0\hspace{-1mm}\bmat{A_{di}(s) & \hspace{-1mm}B_{1di}(s) &\hspace{-1mm}B_{2di}(s)& \hspace{-1mm}E_{di}(s)\\C_{1di}(s) & \hspace{-1mm}D_{11di}(s) & \hspace{-1mm}D_{12di}(s)& \hspace{-1mm}E_{1di}(s)\\C_{2di}(s) &\hspace{-1mm}D_{21di}(s) & \hspace{-1mm}D_{22di}(s)& \hspace{-1mm}E_{2di}(s)} \bmat{x(t+s)\\w(t+s)\\u(t+s)\\ \dot x(t+s)}\hspace{-1mm}ds\label{eqn:NDS}\\[-8mm]\notag
\end{align}
\textbf{The Class of ODE-PDE Systems:}\vspace{-3mm}
\begin{align}
	\bmat{\dot{x}(t)\\ z(t)\\y(t)\\ \phi_i(t,0)}&=\bmat{A_0 & B_1& B_2\\C_1 &D_{11}&D_{12}\\C_2&D_{21}&D_{22}\\C_{ri}&B_{r1i}&B_{r2i}}\bmat{x(t)\\w(t)\\u(t)}+\bmat{B_v\\D_{1v}\\D_{2v}\\D_{rvi}} v(t) \notag\\
\dot \phi_i(t,s)&=\frac{1}{\tau_i}\phi_{i,s}(t,s),\qquad v(t)=\sum_{i=1}^K C_{vi} \phi_i(t,-1)+\sum_{i=1}^K \int_{-1}^0\hspace{-1.5mm}\tau_iC_{vdi}(\tau_i s) \phi_i(t, s)ds\label{eqn:ODEPDE}\\[-8mm]\notag
\end{align}

\textbf{The Class of Partial Integral Equation (PIE) Systems:}\vspace{-3mm}
\begin{align}
\mcl T \dot{\mbf x}(t)+\mcl B_{T_1}\dot w(t)+\mcl B_{T_2}\dot u(t)&=\mcl A\mbf x(t)+\mcl B_1w(t)+\mcl B_2u(t)\notag\\
\hspace{-1cm}z(t)=\mcl C_1\mbf x(t)+\mcl D_{11}w(t)&+\mcl D_{12}u(t),\notag\\
\hspace{-1cm}y(t)=\mcl C_2\mbf x(t)+\mcl D_{21}w(t)&+\mcl D_{22}u(t)\label{eqn:PIE}\\[-7mm]\notag
\end{align}

\caption{Formulation of the DDE, DDF, NDS, ODE-PDE, and PIE Representations of Systems with Delay}\label{fig:representation}.

\vspace{4mm}
\end{figure*}

Finally, we emphasize that this paper does not advocate for any particular time-domain representation (we do not consider the literature on analysis and control in the frequency domain), be it the DDE, DDF, ODE-PDE, or PIE formulation, and does not propose any new algorithms for analysis and control of delay systems per se. Rather, the purpose of this document is to serve as a guide to representation of delay systems in each framework. Specifically, for each representation, we: state the most general form of each representation - allowing for delays in input, output, process and state; define a notion of solution in each case; provide formulae for conversion between representations under which solutions are equivalent; and briefly list advantages and limitations of the representation as applied to network models of the form of Eqn.~\eqref{eqn:network1}. As discussed in the conclusions, these results can be used to establish notions of stability which are equivalent in all representations and to allow for conversion of optimal controllers and estimators between representations.

While subsets of the DDF and ODE-PDE representations of delay systems can be found in the literature~\cite{bensoussan_book,gu_2010,pepe_2008,mazenc_2013,pepe_2008b,gu_2003,niculescu_book}, and some of these equivalences are known~\cite{karafyllis_2014,richard_2003}, previous works do not:  consider all input-output signals and sources of delay; include PIEs; compare the relative advantages of the models as applied to networks; or provide formulae for conversion between representations. This guide, then, may be used as a convenient source of information for researchers interested in either selection of a representation or conversion of a representation to an alternative format. For convenience and comparison, all representations are listed in Figure~\ref{fig:representation}. All conversion formulae are listed in Figures~\ref{fig:formulae1} and~\ref{fig:formulae2}. Finally, note that all proofs have been omitted, but are included in the extended version of this paper on Arxiv~\cite{peet_2020arxiv_TDS}.\vspace{-4mm}

\paragraph*{Notation} $I_n$ is the identity matrix in $\R^{n\times n}$, $e_i$ is the $i^{th}$ canonical unit vector, $\mbf 1_n$ is the dimension $n$ vector of all ones. $0_{n,m}$ is the zero matrix of dimension $\R^{n\times m}$ and $W^{n,2}[X]$ is the nth-order Sobolev subspace of $L_2$[X].\vspace{-4mm}

%%%%%%%%%%%%%%%%%%%%%%%%%%%%%%%%%%%%%%%%%%%%%%%%%%%%%%%%%%%%%%%%%%%%%%%%%%%%%%%%
\section{The DDE Representation}\label{sec:DDE}\vspace{-3mm}
We begin by defining the signals in the Delay-Differential Equation (DDE) representation:\vspace{-3mm}
\begin{itemize}
\item The present state $x(t)\in \R^n$
\item The disturbance or exogenous input, $w(t)\in \R^m$
\item The controlled input, $u(t)\in \R^p$
\item The regulated or external output, $z(t)\in \R^q$
\item The observed or sensed output, $y(t)\in \R^r$\vspace{-2mm}
\end{itemize}
For convenience, we combine all sources of delay (state, input, output, process) into a single set of delays $\{\tau_i\}_{i=1}^K$ with $0<\tau_1<\cdots<\tau_K$. For given $u\in L_2^p$, $w \in L_2^m$, and initial condition $x_0\in W^{1,2}[-\tau_K,0]^n$, we say that $x:[-\tau_K,\infty]\rightarrow \R^n$, $z:[0,\infty]\rightarrow \R^q$, and $y:[0,\infty]\rightarrow \R^r$ satisfy the DDE defined by $\{A_{i}, B_i,C_i,D_{ij},\cdots\}$ if $x$ is differentiable on $[0,\infty]$ (from the right at $t=0$), $x(s)=x_0(s)$ for $s \in [-\tau_K,0]$, and Eqns.~\eqref{eqn:DDE} are satisfied for all $t \ge 0$. If any $B_{1i},D_{11i},D_{21i}\neq0$, we require $w\in W^{1,2}[0,\infty]^m$ and $w(s)=0$ for $s \le 0$. If any $B_{2i},D_{12i},D_{22i}\neq0$, we require $u\in W^{1,2}[0,\infty]^p$ and $u(s)=0$ for $s \le 0$.  %These constraints ensure existence of $w(t-\tau_i)$ and $u(t-\tau_i)$, respectively.

Under the conditions stated above, existence of a classical continuously differentiable solution $x$ is guaranteed as in, e.g. Thm. 3.3 of Chapter 3 in~\cite{kolmanovskii_book} (See also Thm.~1.1 of Chapter 6 in~\cite{hale_book}).
%Note that there are several formulations of the Cauchy problem for this class of systems. For our purposes, we use the following.
%\begin{defn}
%For a given $x_0 \in \mathcal{C}[-\tau_K,0]$, $w\in L_2^m$, and $u \in L_2^p$, we say that $x,y,z$ solve the DDE defined by Eqns.~\eqref{eqn:DDE1}-\eqref{eqn:DDE3} on $[0,T]$ if $x(s)=x_0(s)$ for $s \in [-\tau_K,0]$, Eqns.~\eqref{eqn:DDE1}-\eqref{eqn:DDE3} hold for $t\in [0,T]$ where $\dot x(t)$ is the right-hand Dini derivative.
%\end{defn}
%Establishing existence of solutions of the DDE representation is straightforward using Gronwall-Bellman and the method of steps.
%\begin{lem}
%For any given $x_0 \in L_2[-\tau_K,0]$, $w\in L_2^m$, $u \in L_2^p$ and $T\ge 0$ there exist $x,y,z$ which solve the DDE defined by Eqns.~\eqref{eqn:DDE1}-\eqref{eqn:DDE3} on $[0,T]$.
%\end{lem}
Note that the dimensions of all matrices in this representation can be inferred from the dimension of the respective state and signals.\vspace{-4mm}

\subsection{Advantages of the DDE Formulation}\vspace{-3mm}
The DDE formulation is the prima facie modeling tool for systems with delay and as such is used in almost all network models. The DDE representation has a clear and intuitive meaning. Furthermore, most algorithms and analysis tools are built for this representation. Specifically, Lyapunov-Krasovskii and Lyapunov-Razumikhin stability tests are naturally formulated in this framework.
However, the DDE does not allow for the representation of difference equations and does not allow us to identify which of the states and inputs are delayed by which amount. For this reason, we consider next the DDF representation.\vspace{-4mm}

\section{The DDF Representation}\label{sec:DDF}\vspace{-3mm}
A generalization of the DDE is the Differential-Difference (DDF) formulation. In addition to the signals included in the DDE, the DDF adds the following.
\begin{itemize}
\item The items stored in the signal $r_i(t)\in \R^{p_i}$ are the parts of $x$, $w$, $u$, $v$ which are delayed by amount $\tau_i$. The $r_i$ are the infinite-dimensional part of the system.
\item The ``output'' signal $v(t)\in \R^{n_v}$ extracts information from the infinite-dimensional signals $r_i$ and distributes this information to the state, sensed output, and regulated output. This information can also be re-delayed by feeding back directly into the $r_i$.
\end{itemize}
The governing equations may now be represented in the more compact form of Eqns.~\eqref{eqn:DDF}.
%\begin{align}
%	\dot{x}(t)&=A_0x(t)+B_1w(t)+B_2u(t)+B_v v(t) \label{eqn:DDF}\\
%	z(t)&=C_{1}x(t)+D_{11}w(t)+D_{12}u(t)+D_{1v} v(t)\notag\\
%	y(t)&=C_{2}x(t)+D_{21}w(t)+D_{22}u(t)+D_{2v} v(t)\notag\\
%r_i(t)&=C_{ri}x(t)+B_{r1i}w(t)+B_{r2i}u(t)+D_{rvi}v(t)\notag\\[-1mm]
%v(t)&=\sum_{i=1}^K C_{vi} r_i(t-\tau_i)+\sum_{i=1}^K \int_{-\tau_i}^0C_{vdi}(s) r_i(t+s)ds.\notag\\[-6mm]\notag
%\end{align}

For given $u\in L_2^p$, $w \in L_2^m$, and initial conditions $x_0\in \R^n$, $r_{i0} \in W^{1,2}[-\tau_i,0]^{p_i}$ satisfying the ``sewing condition''\vspace{-2mm}
\begin{align*}
&r_{i0}(0)=C_{ri}x_0\\
&+D_{rvi}\left(\sum_{i=1}^K C_{vi} r_{i0}(-\tau_i)+\sum_{i=1}^K \int_{-\tau_i}^0C_{vdi}(s) r_{i0}(s)ds\right)\\[-7mm]
\end{align*}
for $i=1,\cdots,K$, we say that $x:[0,\infty]\rightarrow \R^n$, $z:[0,\infty]\rightarrow \R^q$, $y:[0,\infty]\rightarrow \R^r$, $r_i:[-\tau_i,\infty]\rightarrow \R^{p_i}$ for $i=1,\cdots,K$, and $v:[0,\infty]\rightarrow \R^{n_v}$ satisfy the DDF defined by $\{A_{i}, B_i,C_i,D_{ij},\cdots\}$ if $x$ is differentiable on $[0,\infty]$, $r_i(s)=r_{i0}(s)$ for $s \in [-\tau_i,0]$, $r_i(t+\cdot) \in W^{1,2}[-\tau_i,0]$ for $i=1,\cdots,K$, and Eqns.~\eqref{eqn:DDF} are satisfied for all $t \ge 0$. In this manuscript, we assume the $C_{vdi}$ are bounded and in the case where $B_{r1i}\neq 0$ or $B_{r2i}\neq 0$, we require $w\in W^{1,2}[0,\infty]^m$ with $w(s)=0$ for $s\le0$ or $u\in W^{1,2}[0,\infty]^p$ with $u(s)=0$ for $s\le0$, respectively.

Under the conditions stated above, existence of a classical solution $x,r_i,v$ is guaranteed as in~\cite{hale_book}, Chapter 9, Thm.~1.1. Furthermore, the ``sewing condition'' and constraints on $w$ and $u$ ensure the solution $r_i$ is continuously differentiable as in~\cite{gil_book} p. 226; or~\cite{kolmanovskii_book}, Thms. 3.1 and 5.4. Note also that the condition $r_i(t+\cdot) \in W^{1,2}$ may be relaxed to continuity as treated in~\cite{henry_1974}.\vspace{-4mm}
%
%The Cauchy problem for this class of systems is formulated so as to be consistent with the DDE class.
%\begin{defn}
%For a given $r_{i0} \in \mathcal{C}[-\tau_i,0]$ for $i=1,\cdots,K$, $w\in L_2^m$, and $u \in L_2^p$, we say that $x,y,z,v,r_i$ solve the DDF defined by Eqns.~\eqref{eqn:DDF} on $[0,T]$ if $r_i(s)=r_{i0}(s)$ for $s \in [-\tau_i,0]$ for $i=1,\cdots,K$, Eqns.~\eqref{eqn:DDF} hold for $t\in [0,T]$ where $\dot x(t)$ is the right-hand Dini derivative.
%\end{defn}
%As for DDEs, existence of solutions of the DDE representation is straightforward using Gronwall-Bellman and the method of steps.
%\begin{lem}
%For any given $r_{i0} \in \mathcal{C}[-\tau_i,0]$ for $i=1,\cdots,K$, $w\in L_2^m$, and $u \in L_2^p$ and $T\ge 0$ there exist $x,y,z,v,r_i$ which solve the DDE defined by Eqns.~\eqref{eqn:DDE1}-\eqref{eqn:DDE3} on $[0,T]$.
%\end{lem}

\subsection{DDEs are a special case of DDFs}\vspace{-3mm}
Although Eqns.~\eqref{eqn:DDF} are more compact, they are more general than the DDEs in~\eqref{eqn:DDE}. Specifically, if we use the conversion formula defined in Eqn.~\eqref{eqn:DDEtoDDF}, then the solution to the DDF is also a solution to the DDE and vice-versa.\vspace{-2mm}
\begin{lem} Suppose that $C_{vi}$, $C_{vdi}$, $C_{ri}$, $B_{r1i}$ , $B_{r1i}$, $D_{rvi}$, $B_v$, $D_{1v}$, and $D_{2v}$ are as defined in Eqns.~\eqref{eqn:DDEtoDDF}. Given $u$, $w$, $x_{0}$, the functions $x$, $y$, and $z$ satisfy the DDE defined by $\{A_{i}, B_i,C_i,D_{ij},\cdots\}$ if and only if $x$, $y$, $z$, and $r_i$ satisfy the DDF defined by $\{A_{i}, B_i,C_i,D_{ij},\cdots\}$ where\vspace{-2mm}
\[
r_i(t)={\bmat{x(t)\\w(t)\\u(t)}},\quad r_{i0}={\bmat{x_0\\0\\0}}\qquad i=1,\cdots,K.\vspace{-6mm}
\]
\end{lem}
%\begin{cor}Suppose $u$, $w$, $x$, $y$, $r_i$ and $z$ satisfy Eqns.~\eqref{eqn:DDF} where $C_{vi}$, $C_{vdi}$, $C_{ri}$, $B_{r1i}$ , $B_{r1i}$, $D_{rvi}$, $B_v$, $D_{1v}$, and $D_{2v}$ are as defined in Eqns.~\eqref{eqn:DDEtoDDF1}-\eqref{eqn:DDEtoDDF2}. Then $u$, $w$, $x$, $y$, and $z$ also satisfy Eqns.~\eqref{eqn:DDE1}-\eqref{eqn:DDE3}.
%\end{cor}

\subsection{Neutral-Delay Systems (NDSs) are a special case of DDFs}\vspace{-3mm}
DDFs are a natural extension of NDSs, which have the general form of Eqn.~\eqref{eqn:NDS} where for simplicity, we assume $x(t),w(t),u(t)=0$ for all $t\le 0$.
%For brevity, we ignore output signals and consider only the following class of neutral systems.\vspace{-2mm}
%\begin{align}
%&\dot x(t)+\sum\nolimits_{i=1}^K E_i \dot x(t-\tau_i)\label{eqn:NDE}\\[-0mm]
%&\quad =Ax(t)+\sum\nolimits_{i=1}^K A_i x(t-\tau_i)+B_1w(t)+B_2u(t)\notag\\[-5mm]\notag
%\end{align}
The conversion from NDS to DDF is given in Eqn.~\eqref{eqn:NDStoDDF}.
\begin{lem} Suppose that $C_{vi}$, $C_{vdi}$, $C_{ri}$, $B_{r1i}$ , $B_{r2i}$, $D_{rvi}$, $B_v$, $D_{1v}$, and $D_{2v}$ are as defined in Eqns.~\eqref{eqn:NDStoDDF}. Given $u$, $w$, the functions $x$, $y$, and $z$ satisfy the NDS defined by $\{A_{i}, B_i,C_i,D_{ij},\cdots\}$ if and only if $x$, $y$, $z$, $v$ and $r_i$ satisfy the DDF defined by $\{A_{i}, B_i,C_i,D_{ij},\cdots\}$ where $r_{i0}=0$ and \vspace{-2mm}
\[
r_i(t)={\bmat{x(t) \\w(t)\\u(t)\\ \dot x(t)}}, \qquad i=1,\cdots,K.\vspace{1mm}
\]
and\vspace{-5mm}
{\small\begin{align*}
&v(t)=\sum_{i=1}^K \bmat{A_i  & B_{1i} & B_{2i}& E_i\\C_{1i}& D_{11i} & D_{12i} & C_{1ei}\\C_{2i} & D_{21i} & D_{22i}&C_{2ei}} \bmat{x(t-\tau_i)\\w(t-\tau_i)\\u(t-\tau_i)\\ \dot x(t-\tau_i)}+\\
&\sum_{i=1}^K \int\limits_{-\tau_i}^0\bmat{A_{di}(s) & B_{1di}(s) &B_{2di}(s)& E_{di}(s)\\C_{1di}(s) & D_{11di}(s) & D_{12di}(s)& C_{1dei}(s)\\C_{2di}(s) &D_{21di}(s) & D_{22di}(s)& C_{2dei}(s)} \bmat{x(t+s)\\w(t+s)\\u(t+s)\\ \dot x(t+s)}ds.
\end{align*}}
\end{lem}

\subsection{Advantages of the DDF Representation}\vspace{-3mm}
The first advantage of the DDF is that it may include difference equations. To illustrate, suppose we set all matrices to zero except $D_{rvi}$ and $C_{vi}$. Then we have the following set of Difference Equations (DEs)\vspace{-2mm}
\[
r_i(t)=\sum\nolimits_{j=1}^K D_{rvi}C_{vj} r_j(t-\tau_j)\qquad i=1,\cdots,K.\vspace{-2mm}
\]
Another example of DEs can be found in Subsec.~\ref{subsec:UAV_SOF}, where we provide a model of network control which can be represented as a DDF, but not a DDE. A related advantage of the DDF is the ability of DDFs to generate discontinuous solutions if the ``sewing condition'' on initial conditions is relaxed. This ability is not inherited using our formulation of ODE-PDE or PIE.

The second advantage of the DDF occurs when the delayed channels only include subsets of the state. For example, if the matrices $A_i$ have low rank (ignoring input and disturbance delay), then $A_{i}=\tilde A_{i}\hat A_i$ for some $\hat A_i$, $\tilde A_i$ where $\hat A_i \in \R^{l_i \times n}$ with $l_i<n$ and we may choose
$C_{vi}=\tilde A_{i}$ and $C_{ri}=A_i$.
 The dimension of $r_i(t)$ now becomes $\R^{l_i}$. This decomposition may be used to reduce complexity in the DDF formulation if $l_i<n$. This reduction is illustrated in detail using the UAV network model in Subsec.~\ref{subsec:UAV_DDF} and the temperature control network in Sec.~\ref{sec:temp}.

A disadvantage of the DDF is that fewer tools are available for analysis and control of DDFs. This is partially because the class of DDFs is larger than the DDEs and thus the tools must be more general. However, we do note that versions of both the Lyapunov-Krasovskii (\cite{gu_2010}) and Lyapunov-Razumikhin (\cite{zhang1998new}) stability tests have been formulated in the DDF framework.\vspace{-3mm}

\section{The Coupled ODE-PDE Representation}\label{sec:ODEPDE}\vspace{-3mm}
We next consider the coupled ODE-PDE representation. Widely recognized as a physical interpretation of delay systems~\cite{richard_2003,hale_book}, ODE-PDE representations allow us to use backstepping methods originally developed for control of PDE models and which have recently been extended to systems with delay - See~\cite{krstic_2008,zhu_2015,karafyllis_2014}.
The particular class of ODE-PDE systems, as given in Eqn.~\eqref{eqn:ODEPDE}, is equivalent to the class of DDFs. Since we have shown that DDEs are a special case of DDFs, we present only the conversion between DDF and ODE-PDE. Such conversion is trivial, however, as all matrices  in the following ODE-PDE model are the same ones used to define the DDF.

For given $u\in L_2^p$, $w \in L_2^m$, and initial conditions $x_0 \in \R^n$, $\phi_{i0} \in W^{1,2}[-1,0]^{p_i}$ satisfying the ``sewing condition''\vspace{-2mm}
\begin{align}
&\phi_{i0}(0)=C_{ri}x_0\label{eqn:sewing_PDE}\\
&+D_{rvi}\left(\sum_{i=1}^K C_{vi} \phi_{i0}(-1)+\sum_{i=1}^K \int_{-1}^0\tau_iC_{vdi}(\tau_i s) \phi_{i0}(s)ds\right)\notag\\[-6mm] \notag
\end{align}
for $i=1,\cdots,K$, we say that $x:[0,\infty]\rightarrow \R^n$, $z:[0,\infty]\rightarrow \R^q$, $y:[0,\infty]\rightarrow \R^r$, $\phi_i(t) \in W^{1,2}[-1,0]^{p_i}$ for $i=1,\cdots,K$, and $v:[0,\infty]\rightarrow \R^{n_v}$ satisfy the ODE-PDE defined by $\{A_{i}, B_i,C_i,D_{ij},\cdots\}$ if $x$ is differentiable and $\phi_i$ is Fr\'echet differentiable on $[0,\infty]$, $x(0)=x_0$, $\phi_i(0,s)=\phi_{i0}(s)$ for $s \in [-1,0]$ for $i=1,\cdots,K$,
and Eqns.~\eqref{eqn:ODEPDE} are satisfied for all $t \ge 0$. As for the DDF, if $B_{r1i}\neq 0$ or $B_{r2i}\neq 0$, we require $w\in W^{1,2}$ or $u\in W^{1,2}$, respectively.\vspace{-2mm}

In Eqns.~\eqref{eqn:ODEPDE}, the infinite-dimensional part of the state is $\phi_i$ - which represents a pipe through which information is flowing. Our formulation is somewhat atypical in that we have scaled all the pipes to have unit length and accelerated or decelerated flow through the pipes according to the desired delay. Solutions to Eqns.~\eqref{eqn:ODEPDE} and Eqns.~\eqref{eqn:DDF} are equivalent, as in the following lemma.\vspace{-2mm}
\begin{lem} Suppose for given $u$, $w$, $r_{i0}$, that $x$, $r_i$, $v$, $y$, and $z$ satisfy the DDF defined by $\{A_{i}, B_i,C_i,D_{ij},\cdots\}$. Then for $u$, $w$, $\phi_{i0}(s)=r_{i0}(\tau_i s)$, we have that $x$, $v$, $y$, and $z$ also satisfy the ODE-PDE defined by $\{A_{i}, B_i,C_i,D_{ij},\cdots\}$ with $
\phi_i(t,s)=r_i(t+\tau_i s)$. Similarly, for given $u$, $w$, $\phi_{i0}$, if $x$, $v$, $y$, $\phi_i$ and $z$ satisfy the ODE-PDE defined by $\{A_{i}, B_i,C_i,D_{ij},\cdots\}$, then $x$, $v$, $y$, and $z$ satisfy the DDF with $r_i(t)=\phi_i(t,0)$ and $r_{i0}(s)=\phi_{i0}(s/\tau_i)$.\vspace{-3mm}
\end{lem}

\subsection{Advantages of the ODE-PDE Representation}\vspace{-3mm}
In the ODE-PDE representation, the infinite-dimensional part of the state is $\mbf \phi(t) \in W^{1,2}[-1,0]^{\sum_i p_i}$. Significantly, by scaling the pipes (and ignoring the distributed delay), the ODE-PDE representation isolates the effect of the delay parameters to a single term - $\dot \phi_i(t,s)=\frac{1}{\tau_i}\phi_{i,s}(t,s)$. This feature makes it easier to understand the effects of uncertainty and time-variation in the delay parameter. Additionally, the ODE-PDE is the native representation used for recently developed backstepping methods for systems with delay, such as proposed in~\cite{krstic_2008,zhu_2015,karafyllis_2014} and use of the conversion formulae provided may allow these methods to be applied to solve a larger class of systems - including difference equations.\vspace{-4mm}

\section{The PIE Representation}\label{sec:PIE}\vspace{-3mm}
A Partial Integral Equation (PIE) has the form of Eqn.~\eqref{eqn:PIE}, where the operators $\mcl T, \mcl A, \mcl B_{i}, \mcl C_i, \mcl D_{ij}$ are Partial Integral (PI) operators and have the form\vspace{-2mm}
\[
		\left(\fourpi{P}{Q_1}{Q_2}{R_i}\bmat{x\\\mbf \Phi}\right)(s):= {\bmat{
				Px + \int_{-1}^{0} Q_1(s)\mbf \Phi(s)ds\\
				Q_2(s)x +\left(\mcl P_{\{R_i\}}\mbf \Phi\right)(s)
		}}\vspace{-2mm}
\]
where\vspace{-2mm}
\begin{align*}
&\left(\mcl{P}_{\{R_i\}}\mbf \Phi\right)(s):= \\
&R_0(s) \mbf \Phi(s) +\int_{-1}^s R_1(s,\theta)\mbf \Phi(\theta)d \theta+\int_s^0R_2(s,\theta)\mbf \Phi(\theta)d \theta.\\[-7mm] \notag
\end{align*}
For given $u\in L_2^p$, $w \in L_2^m$, and initial conditions $\mbf x_{0} \in \R^n \times L_2[-1,0]^{p}$, we say that $\mbf x(t) \in \R^n \times L_2[-1,0]^p$, $z:[0,\infty]\rightarrow \R^q$, $y:[0,\infty]\rightarrow \R^r$ satisfy the PIE defined by $\{\mcl T, \mcl A, \mcl B_i, \mcl C_i, \mcl D_{ij},\mcl B_{T_i}\}$ if $\mbf x$ is Fr\'echet differentiable on $[0,\infty]$, $\mbf x(0)=\mbf x_{0}$ and Eqns.~\eqref{eqn:PIE} are satisfied for all $t \ge 0$. As for the ODE-PDE, if $\mcl B_{T_1}\neq 0$ or $\mcl B_{T_2}\neq 0$ we require $w\in W^{1,2}$ or $u\in W^{1,2}$, with $w(0)=0$ or $u(0)=0$, respectively.\vspace{-2mm}

Heretofore, we have shown that the DDE is a special case of the DDF, which is equivalent to a coupled ODE-PDE, where coupling occurs at the boundary. Given a DDF or ODE-PDE representation, it is relatively straightforward to convert to a PIE by defining the operators  $\mcl T, \mcl A, \mcl B_i, \mcl C_i, \mcl D_{ij},\mcl B_{T_i}$ for which solutions to Eqns.~\eqref{eqn:PIE} also define solutions to Eqns.~\eqref{eqn:DDF} (DDF) and Eqns.~\eqref{eqn:ODEPDE} (ODE-PDE). Specifically, let us define $\{\mcl T, \mcl A, \mcl B_i, \mcl C_i, \mcl D_{ij},\mcl B_{T_i}\}$ as in Eqn.~\eqref{eqn:PIE_ops}
%\vspace{-2mm}
%		\begin{align}
%&\mcl A=\fourpi{\mbf A_0}{\mbf A}{0}{I_\tau,0,0},\quad \mcl T=\fourpi{I}{0}{\mbf T_0}{0,\mbf T_a,\mbf T_b}\label{eqn:PIE_ops}\\
%&\mcl B_1=\fourpi{\mbf B_1}{\emptyset}{0}{\emptyset},\; \mcl B_2=\fourpi{\mbf B_2}{\emptyset}{0}{\emptyset},\notag\\
%& \mcl B_{T_1}=\fourpi{0}{\emptyset}{\mbf T_{1}}{\emptyset}, \;\mcl B_{T_2}=\fourpi{0}{\emptyset}{\mbf T_{2}}{\emptyset},\;\notag\\
%&\mcl C_1=\fourpi{\mbf C_{10}}{\hspace{2mm}\mbf C_{11}}{\emptyset}{\emptyset},\; \mcl C_2=\fourpi{\mbf C_{20}}{\hspace{2mm}\mbf C_{21}}{\emptyset}{\emptyset},\;\mcl D_{ij}=\fourpi{\mbf D_{ij}}{\emptyset}{\emptyset}{\emptyset}\notag\\[-5mm] \notag
%		\end{align}
%		\begin{align}
%\mcl A:=&\fourpi{\mbf A_0}{\mbf A}{0}{I_\tau,0,0},\qquad \mcl T:=\fourpi{I}{0}{\mbf T_0}{0,\mbf T_a,\mbf T_b}\label{eqn:PIE_ops}\\
%\mcl B_1:=&\fourpi{\mbf B_1}{\emptyset}{0}{\emptyset},\qquad \mcl B_2:=\fourpi{\mbf B_2}{\emptyset}{0}{\emptyset}, \\
%\mcl B_{T_1}:=&\fourpi{0}{\emptyset}{\mbf T_{1}}{\emptyset},\qquad
%\mcl B_{T_2}:=\fourpi{0}{\emptyset}{\mbf T_{2}}{\emptyset}\\
%\mcl C_1:=&\fourpi{\mbf C_{10}}{\hspace{2mm}\mbf C_{11}}{\emptyset}{\emptyset},\quad \mcl C_2:=\fourpi{\mbf C_{20}}{\hspace{2mm}\mbf C_{21}}{\emptyset}{\emptyset},\\
%\mcl D_{ij}:=&\fourpi{\mbf D_{ij}}{\emptyset}{\emptyset}{\emptyset}\label{eqn:PIE_opsN}
%		\end{align}
where the required matrices are as defined in Eqns.~\eqref{eqn:PIE_Mats}. Then we have the following.\vspace{-2mm}
\begin{figure*}[!t]
\textbf{Conversion Formula from DDE to DDF:}\vspace{-4mm}
 \begin{equation}
\bmat{B_v\\D_{1v}\\D_{2v}}=I,\;
C_{vi}=\bmat{A_i & B_{1i} &B_{2i}\\C_{1i} & D_{11i}&D_{12i} \\ C_{2i} & D_{21i}&D_{22i}},\;
C_{vdi}(s)=\bmat{A_{di}(s) & B_{1di}(s) &B_{2di}(s)\\C_{1di}(s) & D_{11di}(s)&D_{12di}(s) \\ C_{2di}(s) & D_{21di}(s)&D_{22di}(s)},\;
D_{rvi}=0,\;
\bmat{C_{ri}&B_{r1i}&B_{r2i}}=I.\label{eqn:DDEtoDDF}
\end{equation}
\textbf{Conversion Formula from NDS to DDF:}\vspace{-3mm}
\begin{align}
&D_{rvi}=\bmat{0&0&0\\0&0&0\\0&0&0\\I&0&0},\qquad \bmat{C_{ri}&B_{r1i}&B_{r2i}}=\bmat{I_n & 0& 0\\ 0&I_m&0\\0&0&I_p\\A_0 & B_1 & B_2}  ,\qquad\bmat{B_v\\D_{1v}\\D_{2v}}=I_{n+q+r},\notag\\[-0mm]
&C_{vi}= \bmat{A_i  & B_{1i} & B_{2i}& E_i\\C_{1i}& D_{11i} & D_{12i} & E_{1i}\\C_{2i} & D_{21i} & D_{22i}&E_{2i}},\quad C_{vdi}(s)=\bmat{A_{di}(s) & B_{1di}(s) &B_{2di}(s)& E_{di}(s)\\C_{1di}(s) & D_{11di}(s) & D_{12di}(s)& E_{1di}(s)\\C_{2di}(s) &D_{21di}(s) & D_{22di}(s)& E_{2di}(s)}\label{eqn:NDStoDDF}\\[-6mm]\notag
\end{align}
\textbf{Conversion Formula from ODE-PDE or DDF to PIE:}\vspace{-3mm}
		\begin{align}
\mcl A&=\fourpi{\mbf A_0}{\mbf A}{0}{I_\tau,0,0},& \mcl T&=\fourpi{I}{0}{\mbf T_0}{0,\mbf T_a,\mbf T_b}, & \mcl B_{T_1}&=\fourpi{0}{\emptyset}{\mbf T_{1}}{\emptyset}, &\mcl B_{T_2}&=\fourpi{0}{\emptyset}{\mbf T_{2}}{\emptyset},&\notag\\
\mcl B_1&=\fourpi{\mbf B_1}{\emptyset}{0}{\emptyset}, &\mcl B_2&=\fourpi{\mbf B_2}{\emptyset}{0}{\emptyset},& \mcl C_1&=\fourpi{\mbf C_{10}}{\hspace{2mm}\mbf C_{11}}{\emptyset}{\emptyset},& \mcl C_2&=\fourpi{\mbf C_{20}}{\hspace{2mm}\mbf C_{21}}{\emptyset}{\emptyset},&\mcl D_{ij}=\fourpi{\mbf D_{ij}}{\emptyset}{\emptyset}{\emptyset}\label{eqn:PIE_ops}\\[-7mm] \notag
		\end{align}
where\vspace{-3mm}
\begin{align}
&\hat C_{vi}=C_{vi} +\int_{-1}^0 \tau_iC_{vdi}(\tau_i s)ds,\qquad D_I=\left(I_{n_v}-\left(\sum_{i=1}^K \hat C_{vi}D_{rvi}\right)\right)^{-1},\quad C_{Ii}(s)=-D_I\left(C_{vi} +\tau_i\int_{-1}^sC_{vdi}(\tau_i \eta) d \eta \right)\notag\\
%%%%%%%%%%%%%%%%%%%%%%%%
&\bmat{\mbf T_0 & \mbf T_1 &\mbf T_2}=\bmat{C_{r1} & B_{r11} & B_{r21}\\ \vdots&\vdots&\vdots \\C_{rK} & B_{r1K}& B_{r2K}} +\bmat{D_{rv1} \\ \vdots \\D_{rvK} }\bmat{C_{vx}&D_{vw}&D_{vu}},\quad
\bmat{C_{vx}&D_{vw}&D_{vu}}=D_I\sum_{i=1}^K \hat C_{vi} \bmat{C_{ri}&B_{r1i}&B_{r2i}}\notag\\
%%%%%%%%%%%%%%%%%%%%%%%%
&\mbf T_{a}(s,\theta)=\bmat{D_{rv1} \\ \vdots \\ D_{rvK}}\bmat{C_{I1}(\theta)&\cdots&C_{IK}(\theta)},\qquad \mbf T_{b}(s,\theta)=-I_{\sum_i p_i}+\mbf T_{a}(s,\theta),\quad I_\tau=\bmat{\frac{1}{\tau_1}I_{p_1}&&\\&\ddots&\\&&\frac{1}{\tau_K}I_{p_K}}, \notag\\
 %%%%%%%%%%%%%%%%
&\bmat{\mbf A(s)\\\mbf C_{11}(s) \\\mbf C_{21}(s) }=\bmat{B_v\\D_{1v} \\D_{2v} }\bmat{C_{I1}(s)&\cdots&C_{IK}(s)},\;
\bmat{ \mbf A_0&\mbf B_1&\mbf B_2\\
\mbf C_{10}&\mbf D_{11}&\mbf D_{12}\\
\mbf C_{20}&\mbf D_{21}&\mbf D_{22}}=\bmat{  A_0& B_1&B_2\\
 C_{10}& D_{11}& D_{12}\\
 C_{20}&D_{21}& D_{22}}+\bmat{B_v\\D_{1v} \\D_{2v} }\bmat{C_{vx}&D_{vw}&D_{vu}}.\label{eqn:PIE_Mats}\\[-8mm]\notag
\end{align}
\caption{Conversion formulae from DDE to DDF, NDS to DDF, and DDF/ODE-PDE to PIE}\label{fig:formulae1}
\hrulefill
\end{figure*}

\begin{lem}
  Given $u$, $w$, and $x_0$, $\phi_{i0}$ satisfying the ``Sewing Condition~\eqref{eqn:sewing_PDE}'',  Suppose $x$, $\phi_i$, $v$, $y$, and $z$ satisfy the ODE-PDE defined by $\{A_{i}, B_i,C_i,D_{ij},\cdots\}$. Then $y$ and $z$ also satisfy the PIE defined by $\{\mcl T, \mcl A, \mcl B_i, \mcl C_i, \mcl D_{ij},\mcl B_{T_i}\}$ with $\mcl T, \mcl A, \mcl B_i, \mcl C_i, \mcl D_{ij},\mcl B_{T_i}$ as defined in Eqn.~\eqref{eqn:PIE_ops} and\vspace{-3mm}
\[
\mbf x(t):={\bmat{x(t)\\ \partial_s \phi_{1}(t,\cdot)\\\vdots\\ \partial_s \phi_{K}(t,\cdot)}}\quad \mbf x_0:={\bmat{x_0\\ \partial_s \phi_{10}\\ \vdots \\ \partial_s \phi_{K0}}}.\vspace{-3mm}
\]
Furthermore, for given $u$, $w$, $\mbf x_{0} \in \R^n \times L_2[-1,0]^p$, if $y$, $z$ and $\mbf x$ satisfy the PIE defined by $\{\mcl T, \mcl A, \mcl B_i, \mcl C_i, \mcl D_{ij},\mcl B_{T_i}\}$, then %for
%\[
%\bmat{x_0\\ \phi_{10}\\\vdots\\ \phi_{K0}}=\mcl T \mbf x_0
%\]
 $x$, $\phi_i$, $v$, $y$, and $z$ satisfy the ODE-PDE defined by $\{A_{i}, B_i,C_i,D_{ij},\cdots\}$ where\vspace{-2mm}
\[
{\bmat{x(t)\\ \phi_{1}(t,\cdot)\\\vdots\\ \phi_{K}(t,\cdot)}}=\mcl T \mbf x(t)+\mcl B_{T1}w(t)+\mcl B_{T2}u(t),\; {\bmat{x_0\\ \phi_{10}\\\vdots\\ \phi_{K0}}}=\mcl T \mbf x_0.\vspace{-4mm}
\]
\end{lem}

%\begin{cor}
%  Suppose $u$, $w$, $y$, $\mbf x$ and $z$ satisfy Eqns.~\eqref{eqn:PIE} with $\mcl T, \mcl A, \mcl B_{i}, \mcl C_i, \mcl D_{ij}$ as defined in~\eqref{eqn:PIE_ops}-\eqref{eqn:PIE_opsN}. Then $u$, $w$, $y$, and $z$ satisfy Eqns.~\eqref{eqn:ODEPDE} with
%\[
%\bmat{x(t)\\ \phi_{1}(t,\cdot)\\\vdots\\ \phi_{K}(t,\cdot)}=\mcl T \mbf x(t)+\mcl B_{T1}w(t)+\mcl B_{T2}u(t).
%\]
%%\[
%%\mbf x(t):=\bmat{x(t)\\ \mbf \Phi(t,\cdot)},
%%\]
%\end{cor}
%
%Note that when $D_{rvi}=0$, $D_I=I$.
%\begin{pf}
%See~\cite{peet_2020arxiv_TDS} at Arxiv for the proof.\vspace{-4mm}
%\end{pf}
Note that while solutions of the ODE-PDE are equivalent to those of the PIE, some notions of stability of such solutions may not be. %Specifically, we note that $\mbf x$ and $\phi$$\norm{\mbf x(t)}$
\vspace{-4mm}

\subsection{Advantages of the PIE Representation} \vspace{-3mm}
Like the DDF and ODE-PDE, PIEs can be used to represent low-dimensional delay channels. An additional advantage is the lack of boundary conditions or the `sewing' constraint on the initial condition in, e.g. Eqn.~\eqref{eqn:sewing_PDE}. This is significant in that the implicit dynamics in an ODE-PDE imposed by boundary conditions on $\phi_i$ complicate stability and optimal control problems. By contrast, in PIEs, the infinite-dimension part of the state is $\partial_s \phi_i$ which is in $L_2$ but is otherwise unconstrained. Furthermore, PIEs are defined using the algebra of Partial Integral (PI) operators. The algebraic nature of PI operators implies that most tools developed for matrices can be extended to PIEs - including the LMI framework. Specifically, the LMIs for $H_{\infty}$-optimal observer and controller synthesis have been extended to PIEs, as can be found in~\cite{shuangshuang_2020CDC} and~\cite{shivakumar_2020CDC}, respectively. We refer to Linear PI Inequalities (LPIs) as this extension of the LMI framework and a Matlab toolbox for solving LPIs can be found in~\cite{shivakumar_2020ACC}. An example of these synthesis results can be found in Sec.~\ref{sec:temp}.\vspace{-3mm}

\subsection{Conversion from DDE to PIE}\label{subsec:DDEtoPIE}\vspace{-3mm}
In this subsection, we bypass the DDF and give a formula for direct conversion between the DDE and PIE representations. This formula is given in Eqns.~\eqref{eqn:DDEtoPIE}.\vspace{-4mm}
\begin{figure*}[!t]
\textbf{Conversion Formula from DDE to PIE:} $\mcl T, \mcl A, \mcl B_i, \mcl C_i, \mcl D_{ij},\mcl B_{T_i}$ are as defined in Eqn.~\eqref{eqn:PIE_ops} where now
\begin{align}
&I_\tau=\bmat{\frac{1}{\tau_1}I_{n+m+p}&&\\&\ddots&\\&&\frac{1}{\tau_K}I_{n+m+p}},\; \mbf T_0={\scriptsize \bmat{\bmat{I_{n}&0_{n,m} & 0_{n,p}}^T\\ \vdots \\\bmat{I_{n}&0_{n,m} & 0_{n,p}}^T}},\;
\mbf T_1={\scriptsize \bmat{\bmat{0_{m,n}&I_{m} & 0_{m,p}}^T\\ \vdots \\\bmat{0_{m,n}&I_{m} & 0_{m,p}}^T}},\;  \mbf T_2={\scriptsize \bmat{\bmat{0_{p,n}&0_{p,m} & I_{p}}^T\\ \vdots \\\bmat{0_{p,n}&0_{p,m} & I_{p}}^T}},\notag\\
%%%%%%%%%%%%%%%%%%%%%%%%
&\mbf T_a=0_{(n+m+p)K},\quad \mbf T_b=-I_{(n+m+p)K},\notag \\
%%%%%%%%%%%%%%%
&\bmat{\mbf A(s) \\ \mbf C_{11}(s) \\ \mbf C_{21}(s)}=-\bmat{X_{1}(s)&\cdots&X_{K}(s) },\; \qquad X_{i}(s)=\bmat{A_i & B_{1i} &B_{2i} \\C_{1i} & D_{11i}&D_{12i} \\ C_{2i} & D_{21i}&D_{22i}} +\tau_i\int_{-1}^s\bmat{A_{di}(\tau_i\eta) & B_{1di}(\tau_i\eta) &B_{2di}(\tau_i\eta) \\ C_{1di}(\tau_i\eta) & D_{11di}(\tau_i\eta)&D_{12di}(\tau_i\eta)\\C_{2di}(\tau_i\eta) & D_{21di}(\tau_i\eta)&D_{22di}(\tau_i\eta) } d \eta,\notag\\
%%%%%%%%%%%%%%%%
&\bmat{ \mbf A_0&\mbf B_1&\mbf B_2\\
\mbf C_{10}&\mbf D_{11}&\mbf D_{12}\\
\mbf C_{20}&\mbf D_{21}&\mbf D_{22}}=\bmat{  A_0& B_1&B_2\\
 C_{10}& D_{11}& D_{12}\\
 C_{20}&D_{21}& D_{22}}+\sum_{i=1}^K\bmat{  A_{i}& B_{1i}&B_{2i}\\
 C_{1i}& D_{11i}& D_{12i}\\
 C_{2i}&D_{21i}& D_{22i}}+\int_{-1}^0\sum_{i=1}^K \tau_i \bmat{
 A_{di}(\tau_is)& B_{1di}(\tau_is)&B_{2di}(\tau_is)\\
 C_{1di}(\tau_is)& D_{11di}(\tau_is)& D_{12di}(\tau_is)\\
 C_{2di}(\tau_is)&D_{21di}(\tau_is)& D_{22di}(\tau_is)}ds\label{eqn:DDEtoPIE}\\[-8mm]\notag
\end{align}
\caption{Direct conversion formula from DDE to PIE, bypassing the DDF.}\label{fig:formulae2}
\hrulefill
\end{figure*}

\section{Modeling of a Network of UAVs}\label{sec:UAV} \vspace{-3mm}
To compare the DDE, DDF, ODE-PDE and PIE representations, we return to control of a network of UAVs. In this section, we focus on the DDE and DDF representations, as conversion from DDF to ODE-PDE or PIE is straightforward using the formulae provided. For simplicity, we eliminate the state delays $\hat \tau_{ij}$ governing interactions between UAVs (we will consider state delays in Sec.~\ref{sec:temp}) and map the process, input, and output delays to a common set of delays, $\{\tau_j\}_{j=1}^{3N}$ where the index for the process delay for UAV $i$ is as $\tau_i=\bar \tau_i$, the index for input delay for UAV $i$ is as $\tau_{N+i}=h_i$, and the index of the output delay from UAV $i$ is as $\tau_{2N+i}=\tilde \tau_i$. The process noise is dimension $w(t)\in\R^m$, the common input is dimension $u(t)\in\R^p$, all states are dimension $x_i(t) \in \R^n$ and the outputs are all dimension $y_i(t)\in\R^r$. In this case, we re-write the network model in Eqns.~\eqref{eqn:network1} as\vspace{-2mm}
\begin{align*}
\dot x_i(t)&=a_ix_i(t)+\sum\nolimits_{j=1}^N a_{ij} x_j(t)\\
&\hspace{2cm}+b_{1i}w(t- \tau_{i})+b_{2i} u(t-\tau_{N+i})\\
z(t)&=C_{1} x(t)+D_{12}u(t)\\
y_i(t)&=c_{2i} x_i(t-\tau_{2N+i})+d_{21i}w(t-\tau_{2N+i}).\\[-8mm]
\end{align*}
\subsection{The DDE Representation} \label{subsec:UAV_DDE}\vspace{-3mm} To model this network as a DDE, we consider Eqn.~\eqref{eqn:DDE} where $K=3N$ for a given $C_{10}$ and $D_{12}$.
First, we define $A_0$ blockwise as\vspace{-3mm}
\[
[A_0]_{ij}=\begin{cases}
           a_i, & i=j\\
           a_{ij} & \text{otherwise}
         \end{cases}\vspace{-3mm}
\]
and define the following matrices blockwise for $i=1,\cdots,N$ as\vspace{-3mm}
\begin{align*}
&B_{1,i}= e_i \otimes b_{1i},\quad B_{2,N+i}=e_i \otimes b_{2i},\\
&C_{2,2N+i}= e_i \otimes c_{2i},\quad D_{21,2N+i}=e_i \otimes d_{2i}.\\[-6mm]
\end{align*}
All other undefined matrices in Eqn.~\eqref{eqn:DDE} are $0$.
%Thus, for example, $B_{1,(N+1)}=0$, $B_{2,N}=0$ and $C_{2,2N}=0$.
%\[
%B_{11}=\bmat{b_{11}\\0\\ \vdots \\ 0}, B_{12}=\bmat{0\\b_{12}\\0\\ \vdots \\ 0}, B_{1N}=\bmat{0\\ \vdots \\b_{1N} }
%\]
%\[
%B_{2,N+1}=\bmat{b_{2,1}\\ \vdots \\ 0}, B_{2,2N}=\bmat{0\\  \vdots \\ b_{2,N}},
%\]
%\[
%C_{2,2N+1}=\bmat{c_{21}\\ 0 \\ \vdots \\ 0}, D_{21,2N+1}=\bmat{d_{21}\\ 0 \\ \vdots \\ 0}\text{et c.}
%\]
The DDE representation of the network has the obvious disadvantage that there are $3N$ delays and each delayed channel contains all states and inputs - yielding an aggregate delayed channel of size $\R^{3N(nN+m+p)}$.\vspace{-3mm}

\subsection{The DDF Representation} \label{subsec:UAV_DDF}\vspace{-3mm} To efficiently model the network model as a DDF, we retain the matrix $A_0$ from the DDE model in Subsec.~\ref{subsec:UAV_DDE}, set $C_1=C_{10}$ and leave $D_{12}$ unchanged.
Our first step is to define the vectors $r_i(t)$ and $v(t)$ using $B_{r1i}$, $B_{r2i}$,$C_{ri}$, $C_{vi}$, $B_v$, and $B_{2v}$ (all other matrices are $0$). The first $3$ sets of matrices are defined for $i=1,\cdots,N$  as
%\begin{align*}
%&[B_{r1i}]_{in+1:(i+1)n,jm+1:(j+1)m}=b_{1i} \;  i =N1(j),\;\; j\in [1,N],\\
%&[B_{r2i}]_{in+1:(i+1)n,jp+1:(j+1)p}=b_{2i} \quad  i =N2(j),\;\; j\in [1,N],\\
%&[B_{r1i}]_{in+1:(i+1)n,jp+1:(j+1)p}=b_{1i} \quad  i =N3(j),\;\; j\in [1,N],\\
%&[C_{ri}]_{in+1:(i+1)n,jn+1:(j+1)n}=c_{2i} \quad  i =N3(j),\;\; j\in [1,N].\\
%\end{align*}
$B_{r1,i}= b_{1i}$, $B_{r1,2N+i}= d_{21i}$, $B_{r2,N+i}=b_{2i}$, and $C_{r,2N+i}= c_{2i}$.
%\begin{align*}
%&B_{r1,i}= b_{1i},\qquad B_{r1,2N+i}= d_{21i}, \\
%&B_{r2,N+i}=b_{2i},\qquad C_{r,2N+i}= c_{2i}.
%\end{align*}
We presume the UAV state dimensions ($n$) are less than the size of the aggregate input ($m$) and disturbance vectors ($p$) (i.e. $n<m$ and $n<p$). In this case it is preferable to delay only the part of the input and disturbance signals which affects each UAV. We now have the following definition for $r_i$ for $i=1,\cdots,3N$.\vspace{-2mm}
\begin{align*}
&r_i(t)=\\
&\begin{cases}
b_{1i}w(t)& i\in [1,N]\\
b_{2,i-N}u(t)& i\in [N+1,2N]\\
c_{2,i-2N}x_{i-2N}(t)+d_{21,i-2N}w(t) &  i\in [2N+1,3N].
\end{cases}\\[-5mm]
\end{align*}
Next, we construct output $v(t)$ by defining $C_{vi}$ for $i=1,\cdots,3N$ as
$C_{vi}=e_i\otimes I_{p_i}$
%\[
%[C_{vi}]_j=\begin{cases}I & i=j\\0& \text{otherwise,}\end{cases}\vspace{-2mm}
%\]
%so that
%\[
%C_{v1}=\bmat{I\\0\\\vdots \\ 0}, C_{v2}=\bmat{0\\I\\0\\\vdots \\ 0},\cdots,C_{v,3N}=\bmat{0\\\vdots \\ 0\\I}
%\]
which yields\vspace{-2mm}
\[
v(t)=\bmat{r_1(t-\tau_1)^T& \cdots & r_{3N}(t-\tau_{3N})^T}^T.\vspace{-2mm}
\]
Finally, we feed $v(t)$ back into the dynamics using\vspace{-2mm}
\[
B_v=\bmat{I&\cdots&I&I&\cdots&I&0},\; D_{2v}=\bmat{0&\cdots&0&0&\cdots&0&I},\vspace{-2mm}
\]
which recovers the network model.\vspace{-3mm}
\subsection{Complexity of DDEs vs. DDFs} \vspace{-3mm}
In the DDF model, the infinite-dimensional state is $r_i$. In our DDF formulation of the UAV model: each process delay adds $n$ states; each input delay adds $n$ states; and each output delay adds $r$ states to this vector. The aggregated infinite-dimensional state is then $L_2\text{\textasciicircum}(\sum_{p_i}=(2n+r)N)$. Assuming that optimal control and estimation problems are tractable when the number of infinite-dimensional states is less than 50~\cite{peet_2020SICON}, and if we suppose $n=r=1$, then it is possible to control $17$ UAVs. By contrast, in the DDE model of our UAVs, the infinite-dimensional state is $L_2^{3N(m+p+r)}$ ( meaning we can control at most $5$ or $6$ UAVs).  \vspace{-4mm}%Note that if we had included state delays, the DDF could handle $7$ UAVs, while the DDF could only handle $2$.\vspace{-4mm}
%This number would be further reduced if each UAV has its own input and disturbance (a likely scenario).
%To further illustrate the importance of converting efficiently, we note that if we had included the $N^2$ state delays from the original model, the dimension of $r$ in the DDF formulation would be $(2n+r)N+nN^2$. This may seem large, but if again $n=r=1$, we would still be able to control $6$ or $7$ UAVs (for $N=7$, the dimension is $70$). By contrast, if we had used the naive conversion between DDE and DDF, this dimension would be $(nN+m+p+r)(N^2+3N)$ - meaning we would only be able to control $2$ UAVs using a single input and disturbance (for $N=3$ the dimension is $108$).\vspace{-4mm}

\section{A Network which is a DDF, but not a DDE}\label{subsec:UAV_SOF} \vspace{-3mm}
In this subsection, we present a network model which can be represented using DDFs, ODE-PDEs, and PIEs, but not using DDEs. These models arise from the use of static feedback - i.e. $u(t)=Fy(t)$ where $y(t)$ is the concatenated vector of outputs from the UAVs. Note that $y$ may include measurement of all states (the static state feedback problem). In this example, let us ignore output, process and state delay, but retain input delay and add a term which models the impact of actuator input $u(t)$ on the sensors as\vspace{-2mm}
\[
y_i(t)=c_{2i} x_i(t)+d_{21i}w(t)+d_{22i}u(t-\tau_i).\vspace{-2mm}
\]
Let $A_0$, $C_1$, $D_{12}$, $B_{2i}$, $C_{vi}$ be as defined in Subsec.~\ref{subsec:UAV_DDF} and define\vspace{-2mm}
\[
B_{1}={\bmat{b_{11}\\ \vdots \\ b_{1N}}},\quad D_{21}={\bmat{d_{21,1}\\ \vdots \\ d_{21,N}}}\vspace{-2mm}
\]
\[
C_2=\diag(c_{2,1},\cdots,c_{2,N}),\quad   D_{22i}=e_i \otimes d_{22i}.\vspace{-0mm}
\]
Aggregating the measurements, we have\vspace{-2mm}
\[
y(t)=C_2 x(t)+D_{21}w(t)+\sum\nolimits_{i=1}^N D_{22i}u(t-\tau_i).\vspace{-2mm}
\]
Now, substituting $u(t)=Fy(t)$ into the sensed output term, we obtain solutions of the form\vspace{-2mm}
\begin{align}
\dot x(t)&=A_0x(t)+B_{1}w(t)+\sum\nolimits_{i=1}^N B_{2i} F y(t-\tau_i)\notag \\
z(t)&=C_{1} x(t)+D_{12}Fy(t)\label{eqn:UAV_SOF} \\
y(t)&=C_{2} x(t)+D_{21}w(t)+\sum\nolimits_{i=1}^N D_{22i}Fy(t-\tau_i).\notag\\[-8mm]\notag
\end{align}
%\begin{align*}
%\dot x(t)&=A_0x(t)+B_{1}w(t)+\sum_i B_{2i} u(t-\tau_i)\\
%&=A_0x(t)+B_{1}w(t)+\sum_i B_{2i} F y(t-\tau_i)\\
%z(t)&=C_{1} x(t)+D_{12}u(t)\\
%&=C_{1} x(t)+D_{12}Fy(t)\\
%y(t)&=C_{2} x(t)+D_{21}w(t)+\sum_{i=1}^N D_{22i}u(t-\tau_i)\\
%&=C_{2} x(t)+D_{21}w(t)+\sum_{i=1}^N D_{22i}Fy(t-\tau_i)
%\end{align*}
Clearly, there is no DDE model with solutions which satisfy Eqns.~\eqref{eqn:UAV_SOF} due to the recursion in the output~\cite{henry_1974}. However (assuming appropriate initial conditions), these solutions can be constructed using the DDF (and consequently the ODE-PDE and PIE frameworks). To construct such a model, we define the following terms.\vspace{-2mm}
\begin{align}
\tilde D_{12}&=D_{12}F D_{21},\quad \tilde D_{22}=0,\quad \tilde C_1 =C_1+D_{12}F C_2\notag\\
C_{ri}&=FC_2, \qquad B_{r1i}=FD_{21}, \qquad [D_{rvi}]_i=FD_{22i}\notag\\
B_{v}&=\bmat{B_{21}& \cdots & B_{2N}},\quad C_{vi}=e_i \otimes I
\notag \\
D_{1v}&=D_{12}F D_{2v},\quad D_{2v}=\bmat{D_{22,1}& \cdots & D_{22,N}}\label{eqn:UAV_SOF_defs}\\[-7mm] \notag
\end{align}
\begin{lem}
For given $r_{i0}$, $x_0$, suppose $r_i$, $v$, $y$, $x$, and $z$ satisfy the DDF defined by\vspace{-2mm}
\[
\{A_0,B_1,B_v,\tilde C_{1},\tilde D_{12},D_{1v},D_{2v},C_{ri},B_{r1i},D_{rvi},C_{vi}\}\vspace{-2mm}
\] given by Eqns.~\eqref{eqn:UAV_SOF_defs}. Then $x$, $z$ and $y$ also satisfy Eqns.~\eqref{eqn:UAV_SOF}.
\end{lem}
\vspace{-5mm}
%\begin{pf}
%See~\cite{peet_2020arxiv_TDS} at Arxiv for the proof.\vspace{-2mm}
%\end{pf}
%\begin{pf}
%Suppose $r_i$, $v$, $y$, $x$ satisfy Eqns.~\eqref{eqn:DDF}. Then
%\begin{align*}
%y(t)&=C_{2} x(t)+D_{21}w(t)+D_{2v}v(t)\\
%&=C_{2} x(t)+D_{21}w(t)+\sum\limits_{i=1}^N D_{22i}v_i(t)
%%&=C_{2} x(t)+D_{21}w(t)+\sum_i D_{22i}r(t-\tau_i)
%\end{align*}
%and hence
%\begin{align*}
%r_i(t)&=C_{ri} x(t)+B_{r1i}w(t)+D_{rvi}v(t)\\
%&=FC_{2} x(t)+FD_{21}w(t)+\sum\limits_{i=1}^N FD_{22i}v_i(t)\\
%&=F\left(C_{2} x(t)+D_{21}w(t)+\sum_{i=1}^N D_{22i}v_i(t)\right)=Fy(t).
%\end{align*}
%Next,
%\[
%v_i(t)=r_i(t-\tau_i)=Fy(t-\tau_i)
%\]
%and we conclude
%\[
%y(t)=C_{2} x(t)+D_{21}w(t)+\sum_{i=1}^N D_{22i}Fy(t-\tau_i).
%\]
%Similarly
%\begin{align*}
%\dot x(t)&=A_0x(t)+B_{1}w(t)+B_{v} v(t)\\
%&=A_0x(t)+B_{1}w(t)+\sum\limits_{i=1}^NB_{2i}v_i(t)\\
%&=A_0x(t)+B_{1}w(t)+\sum\limits_{i=1}^NB_{2i}Fy(t-\tau_i)
%\end{align*}
%and finally,
%
%\begin{align*}
%&z(t)=\tilde C_1x(t)+\tilde D_{12}w(t)+D_{1v}v(t)\\
%%&=\left(C_1+D_{12}F C_2\right)x(t)+D_{12}FD_{21}w(t)+\sum_i D_{12}F D_{22i}v_i(t)\\
%&=C_1x(t)+D_{12}F\left(C_{2} x(t)+D_{21}w(t)+\sum\limits_{i=1}^N D_{22i}v_i(t)\right)\\
%&=C_1x(t)+D_{12}Fy(t)
%\end{align*}
%as desired.
%\end{pf}
%
%Note that if $u(t)\in \R^p$, the dimension of the infinite-dimensional state is $\R^{pN}$.

\section{Optimal Control of a Large Network}\label{sec:temp} \vspace{-3mm}
 To illustrate the computational advantages of the DDFs, ODE-PDEs, and PIEs for controller synthesis problems, we consider the scalable network model with state-delay for centralized control of water temperature for multiple showering customers as defined in~\cite{peet_2020SICON}. If $T_{1i}$ is the tap position and $T_{2i}$ is the temperature for user $i$, then the dynamics of this model are given by\vspace{-2mm}
%
%Specifically, we model a user attempting to achieve a desired temperature by adjusting a hot-water tap.  The model assumes user $i$ will adjust their tap position ($T_{1i}(t)$)  at rate proportional to the difference between current temperature ($T_{2i}(t)$) and desired temperature ($w_i(t)$) and with constant of proportionality $\alpha_i$.  When available hot water pressure is finite, the actions of each user affect the temperature of all other users. This is modeled using $\gamma_{ij}$, which represents the fractional reduction of user $i$'s hot water caused by an increase in hot water usage by user $j$. There is also a transport delay, $\tau_i$, caused by flow of hot water from source to the showerhead of user $i$. We add a centralized tracking controller which senses tap position and water temperature. However, the controller can not sense the desired water temperatures, $w_i(t)$ - which is modeled as a disturbance. The regulated output is sum of the tap actions of all users: $T_{2i}$ and the sum of the centralized interventions, $u_i$.
\begin{align}
&\dot T_{1i}(t)=T_{2i}(t)-w_i(t) \label{eqn:shower}\\
&\dot T_{2i}(t)=-\alpha_i \left( T_{2i}(t-\tau_i) -w_i(t)\right) \notag\\
&\qquad + \sum\nolimits_{j \neq i}^N \gamma_{ij} \alpha_{j}\left(T_{2j}(t-\tau_j)-w_j(t)\right)  + u_i(t) \notag\\[-1mm]
 & z(t)  = \bmat{\sum\nolimits_{i=1}^N T_{1i}(t)& .1 \sum\nolimits_{i=1}^N u_i(t)}^T.\notag\\[-6mm]\notag
\end{align}
For $N$ users, we choose $\alpha_i=1$, $\gamma_{ij}=1/N$, $\tau_i=i$, and $w_i(t)=N$.\vspace{-3mm} %We find that for these values, the optimal closed-loop gain from disturbance to regulated output remains relatively constant in the range of $.35-.4$, irrespective of the number of users. \vspace{-2mm}

\subsubsection{DDE Formulation of the Network}\vspace{-3mm}
In~\cite{peet_2020SICON}, we formed the aggregate state vector as\vspace{-2mm}
\[
x(t)=\bmat{T_{11}(t)& \cdots & T_{1N}(t)& T_{21}(t)&\cdots &T_{2N}(t)}^T\vspace{-2mm}
\]
and defined the DDE model using\vspace{-2mm}
\begin{align*}
A_0&=\bmat{0_{N \times N} & I_N\\ 0_{N \times N} & 0_{N \times N}},\quad A_i=\bmat{0_{N \times N} & 0_{N \times N}\\ 0_{N \times N} & \hat A_i}\\
\hat A_i&=\Gamma*\diag(e_i)=\Gamma*\diag\left(\bmat{0_{1 \times i-1} &1& 0_{1 \times N-i}}\right)\\
%\hat A_1&=\bmat{-\alpha_1 &0&0\\vdots &\ddots &0\\ \alpha_1 \gamma_{K1}&0&0}\\
B_1&=\bmat{-I_N\\-\Gamma},\;\qquad B_2=\bmat{0_{N \times N}\\I_N}\\
[\Gamma]_{ij}&=\begin{cases}\gamma_{ij}\alpha_j&i\neq j\\ -\alpha_i & i=j\end{cases} \qquad i,j=1,\cdots,N\\
%(q_i)_j&=\bmat{\alpha_i\gamma_{1i}\\ \vdots \\ \alpha_i\gamma_{Ki}}\\
C_{1}&=\bmat{\mbf 1_N^T&0_{1 \times N}\\0_{1 \times N}&0_{1 \times N}},\;\; D_{11}=\bmat{0_{2 \times N}},\; D_{12}=\bmat{0_{1 \times N}\\.1 \mbf 1_N^T}.\\[-8mm]\notag
\end{align*}
%where $\mbf 1_N$ is the length-$N$ vector of all ones.
In this formulation, we have $n=2N$ states, $m=N$ disturbances, $p=N$ inputs, $q=2$ regulated outputs and $K=N$ delays ($\tau_{ij}=\tau_j$). Using the SOS-based $H_\infty$-optimal controller synthesis algorithm for DDEs as presented in~\cite{peet_2020SICON}, we were able to design controllers for $N=4$ users. This corresponds to an infinite-dimensional channel of size $L_2^{nK=32}$.\vspace{-3mm}

\subsubsection{DDE Formulation of the problem}\vspace{-3mm}
To construct the DDF formulation of the problem, $x(t)$ is unchanged. However, we now define the delayed channels as\vspace{-2mm}
\[
r_i(t)=\bmat{0_{1\times N+i-1}&1&0_{1 \times N-i}}x(t)=T_{2i}(t).\vspace{-2mm}
\]
This is done by defining $C_{ri},B_{r1i}, B_{r2i} $ and $D_{rvi}$ as\vspace{-2mm}
\begin{align*}
C_{ri}&=\bmat{0_{1\times N+i-1}&1&0_{1 \times N-i}}\\
B_{r1i}&=0_{1 \times N}\;\qquad B_{r2i}=0_{1 \times N}\;\qquad D_{rvi}=0_{1 \times N}.\\[-7mm]\notag
\end{align*}
We would like the output of the delayed channels to be the delayed states as\vspace{-2mm}
\[
v(t)=\bmat{T_{21}(t-\tau_1)&\cdots &T_{2N}(t-\tau_N)}^T.\vspace{-2mm}
\]
This is accomplished by defining\vspace{-2mm}
\[
C_{vi}=e_i=\bmat{0_{1 \times i-1}&1&0_{1 \times N-i}}^T,\qquad C_{vdi}=0_{2\times N}.\vspace{-2mm}
\]
Finally, we retain $A_0,B_1,B_2,C_1,C_2,D_{11},D_{12}$ from the DDE formulation, and use $B_v$ and $D_{1v}$  to model how the delayed terms affect the state dynamics and output signal. \vspace{-2mm}
\[
B_v=\bmat{0_{N \times N}\\ \Gamma},\qquad D_{1v}=0\vspace{-2mm}
\]
In the DDF formulation, we have $n=2N$ states, $m=N$ disturbances, $p=N$ inputs, $q=2$ regulated outputs,  $K=N$ delays ($\tau_{ij}=\tau_j$), and $K$ delay channels, each of dimension $L_2^1$.\vspace{-2mm}

\subsection{$H_\infty$-optimal Control Using PIETOOLS 2020a} \vspace{-2mm}
For $H_\infty$-optimal controller synthesis, we used the DDF to PIE converter \texttt{convert\_PIETOOLS\_DDF} and $H_\infty$-optimal synthesis option in the PIETOOLS 2020a Matlab toolbox, as described in~\cite{shivakumar_2020ACC} and available online at~\cite{PIETOOLS_website}. The DDF system input format for this toolbox is described in detail in the user manual~\cite{PIETOOLS_website}, as is the converter and controller synthesis feature. In this toolbox, the extreme performance option was selected to decrease computation times and reduce memory usage. The $H_\infty$-optimal controller synthesis feature in PIETOOLS solves the optimal control problem for a PIE and is based on the result in~\cite{shivakumar_2020CDC}. The numerical test was performed on a desktop computer with 128GB RAM and a 3 GHz intel processor. CPU seconds is as listed for the interior-point calculations determined by Sedumi.  The computation times, indexed by number of users, are listed in Table~\ref{tab:computation}. In all cases, the achieved closed-loop $H_\infty$-norm was in the interval $[.3,3]$. Practically, we observe that the controller synthesis problem is tractable up to 40 users - a significant improvement from the $4$ users in~\cite{peet_2020SICON}. Note that 40 users corresponds to an aggregated infinite-dimensional channel of size $L_2^{\sum_i p_i=N=40}$. Also recall that for 40 users, we have 80 states, 40 inputs, 40 disturbances and 40 delays.

Note that the PIETOOLS 2020a toolbox does not require use of the PIE formulation and will convert a DDE to a DDF, if desired.  There is also a feature for constructing minimal DDF representations of DDEs - which can be very useful for solving large network problems. The conversion from a NDS to DDF is also included in the PIETOOLS library \texttt{examples\_DDF\_library\_PIETOOLS.m}.

%\begin{table}
%\caption{CPU sec indexed by \# of states ($N$)}\label{tab:computation}
%{\scriptsize
%\parbox{.99\linewidth}{
%  \centering
%\begin{tabular}{c|c|c|c|c|c|c|c}
%%&\multicolumn{2}{c}{SOS~\cite{peet_2009SICON}} & \multicolumn{2}{c}{SOS-joint}\\
%%&\multicolumn{2}{c}{SOS~\cite{peet_2009SICON}} & \multicolumn{2}{c}{SOS-joint}\\
%{\hspace{-2mm}\tiny \text{\# of users} }\hspace{-2mm} & $1$  & $3$ & $5$ & $10$ & $20$ & $30$ & $40$\\
%\hline
%\tiny \text{IPM CPU sec} & .48 & .638 & 2.42  &  94.7 & 5455 & 35k&157k\\
%\end{tabular}\vspace{2mm}
%}
%}
%\end{table}
\begin{table}
\caption{IPM CPU sec vs. \# of states ($N$) for $H_\infty$ control of Eqn.~\eqref{eqn:shower}.}
\label{tab:computation}
\begin{center}{\scriptsize
\begin{tabular}{c|c|c|c|c|c|c|c}
%&\multicolumn{2}{c}{SOS~\cite{peet_2009SICON}} & \multicolumn{2}{c}{SOS-joint}\\
%&\multicolumn{2}{c}{SOS~\cite{peet_2009SICON}} & \multicolumn{2}{c}{SOS-joint}\\
{\hspace{-2mm}\tiny \text{N}$\rightarrow$ }\hspace{-2mm} & $1$  & $3$ & $5$ & $10$ & $20$ & $30$ & $40$\\
\hline
\tiny \text{CPU sec} & .48 & .638 & 2.42  &  94.7 & 5455 & 35k&157k\\
\end{tabular}}
\end{center}
\end{table}\vspace{-3mm}

\section{Conclusion}\vspace{-3mm}  This paper summarizes four possible representations for systems with delay: the Delay-Differential Equation (DDE) form; The Differential Difference (DDF) form; the ODE-PDE form; and the Partial Integral Equation (PIE) form. Formulae are given for conversion between these representations, although direct conversion from DDE to DDF is not advised if the delayed channels are low-dimensional (although PIETOOLS 2020a includes a feature for constructing minimal DDF representations of DDEs). Using the given formulae and definitions of solution, we show that the set of solutions for the DDF, ODE-PDE, and PIE are equivalent. These results imply that if there is a valid conversion formula, many solutions to the $H_\infty$-optimal control and estimation problems can be converted between representations by applying this formula to the closed-loop system. However, this only works if optimality is defined in terms of the finite-dimensional vectors, $x_0,u,w,x,y,z$. This is because any input-output pair $(u,w,x_0) \mapsto (y,z,x)$ which defines a solution to one representation also defines a solution for every other representation for which there is a valid conversion formula. Likewise, stability of the representations is equivalent as long as the stability definition only involves the finite-dimensional vectors, $x_0,x,u,w,y,z$.

The results and formulae in this paper are meant to provide a convenient reference for researchers interested in exploring alternative representations of delay systems. A summary of the representations and conversion formulae is given in Table~\ref{tab:summary}, along with examples of simulation tools and controller synthesis results. We have shown using an example of a network of UAVs that some networks cannot be modeled in the DDE formulation and that careful choice of representation can significantly reduce the complexity of the underlying analysis and control problems. Finally, we have shown that $H_\infty$-optimal control in the DDF/ODE-PDE/PIE framework allows up to 40 agents, while formulation in the DDE framework only allows for control of 4 agents.
\begin{table}[hb]
\caption{Conversion formulae (DDF,PDE,PIE), simulation tools (Sim), controller design tools ($H_\infty$), and model definitions (Model) for each class of systems (PDE$\rightarrow$ODE-PDE).}
\label{tab:summary}
\begin{center}{\scriptsize
\begin{tabular}{c|c|c|c|c|c|c}
%&\multicolumn{2}{c}{SOS~\cite{peet_2009SICON}} & \multicolumn{2}{c}{SOS-joint}\\
%&\multicolumn{2}{c}{SOS~\cite{peet_2009SICON}} & \multicolumn{2}{c}{SOS-joint}\\
{\hspace{-2mm}\tiny Need$\rightarrow$ }\hspace{-2mm} & DDF  & PDE & PIE & Sim & $H_\infty$ &Model\\
\hline
\tiny DDE & \eqref{eqn:DDEtoDDF} & \eqref{eqn:DDEtoDDF}+\eqref{eqn:ODEPDE} & \eqref{eqn:DDEtoPIE}  &  \cite{bellen_book}& \cite{peet_2020SICON} & \eqref{eqn:DDE}\\
\tiny Neut. & \eqref{eqn:NDStoDDF} & \eqref{eqn:NDStoDDF}+\eqref{eqn:ODEPDE} & \eqref{eqn:NDStoDDF}+\eqref{eqn:PIE_Mats}  &  \cite{bellen_book} & \cite{xu_2003} & \eqref{eqn:NDS}\\
\tiny DDF & X & \eqref{eqn:ODEPDE} & \eqref{eqn:PIE_Mats}  &  - & - & \eqref{eqn:DDF}\\
\tiny PDE & X & X & \cite{shivakumar_2019CDC}  &  \cite{wouwer_book} & \cite{krstic_2008} & \eqref{eqn:ODEPDE}\\
\tiny PIE & X & X & X  &  - & \cite{shivakumar_2020CDC} & \eqref{eqn:PIE}\\
\end{tabular}}
\end{center}\vspace{-2mm}
\end{table}
\vspace{-4mm}

%\begin{table}[hb]
%\caption{The caption comes before the table.}
%\begin{center}
%\begin{tabular}{|c||c|c|}\hline
%$\otimes$&0&1\\\hline\hline
%0&0&1\\\hline
%1&1&0\\\hline
%\end{tabular}
%\end{center}
%\end{table}

\bibliographystyle{plain}
\bibliography{peet_bib,delay}
\newpage
\newpage

\appendix{Proof of \textbf{Lemma 1}}
% \begin{equation}
%\bmat{B_v\\D_{1v}\\D_{2v}}=I,\;
%C_{vi}=\bmat{A_i & B_{1i} &B_{2i}\\C_{1i} & D_{11i}&D_{12i} \\ C_{2i} & D_{21i}&D_{22i}},\;
%C_{vdi}(s)=\bmat{A_{di}(s) & B_{1di}(s) &B_{2di}(s)\\C_{1di}(s) & D_{11di}(s)&D_{12di}(s) \\ C_{2di}(s) & D_{21di}(s)&D_{22di}(s)},\;
%D_{rvi}=0,\;
%\bmat{C_{ri}&B_{r1i}&B_{r2i}}=I.\label{eqn:DDEtoDDF}
%\end{equation}

\begin{lem}[Lemma 1] Suppose that $C_{vi}$, $C_{vdi}$, $C_{ri}$, $B_{r1i}$ , $B_{r1i}$, $D_{rvi}$, $B_v$, $D_{1v}$, and $D_{2v}$ are as defined in Eqns.~\eqref{eqn:DDEtoDDF}. Given $u$, $w$, $x_{0}$, the functions $x$, $y$, and $z$ satisfy the DDE defined by $\{A_{i}, B_i,C_i,D_{ij},\cdots\}$ if and only if $x$, $y$, $z$, and $r_i$ satisfy the DDF defined by $\{A_{i}, B_i,C_i,D_{ij},\cdots\}$ where
\[
r_i(t)={\bmat{x(t)\\w(t)\\u(t)}},\quad r_{i0}={\bmat{x_0\\0\\0}}\qquad i=1,\cdots,K.\vspace{-4mm}
\]
\end{lem}

\begin{pf} Most of the proof follows immediately from the proof of Lemma 2. However, in this case, we a non-zero initial condition. Specifically, suppose $x$, $y$, and $z$ satisfy the DDE for given $u$, $w$, $x_{0}$. If\vspace{-3mm}
\[
r_i(t)={\bmat{x(t)\\w(t)\\u(t)}},\quad r_{i0}(s)={\bmat{x_0(s)\\0\\0}},\vspace{-3mm}
\]
then as in the proof of Lemma 2, Eqn.~\eqref{eqn:DDF} is satisfied. Furthermore, since $w(t),u(t)=0$ for $t\le 0$, we have for $t\le 0$,\vspace{-3mm}
\[
r_i(s)=\bmat{x(s)\\0\\0}=\bmat{x_0(s)\\0\\0}=r_{i0}(s).\vspace{-3mm}
\]
Finally, $x_{0}\in W^{1,2}$ implies $r_{i0}\in W^{1,2}$ and since $D_{rvi}=0$, the sewing condition is satisfied since\vspace{-3mm}
\[
r_{i0}(0)=C_{ri}x_0=\bmat{x_0\\0\\0}.\vspace{-3mm}
\]
Conversely, suppose $x$, $y$, $z$, and $r_i$ satisfy the DDF. Then from the proof of Lemma 2, Eqn.~\eqref{eqn:DDE} is satisfied. Likewise since $r_{i0} \in W^{1,2}[-\tau_i,0]^{p_i}$ for all $i=1,\cdots,K$, we have $x_0\in W^{1,2}[-\tau_K,0]^n$ which completes the proof.

\end{pf}

\newpage
\appendix{Proof of \textbf{Lemma 2}}

\begin{lem}[Lemma 2] Suppose that $C_{vi}$, $C_{vdi}$, $C_{ri}$, $B_{r1i}$ , $B_{r2i}$, $D_{rvi}$, $B_v$, $D_{1v}$, and $D_{2v}$ are as defined in Eqns.~\eqref{eqn:NDStoDDF}. Given $u$, $w$, the functions $x$, $y$, and $z$ satisfy the NDS defined by $\{A_{i}, B_i,C_i,D_{ij},\cdots\}$ if and only if $x$, $y$, $z$, $v$ and $r_i$ satisfy the DDF defined by $\{A_{i}, B_i,C_i,D_{ij},\cdots\}$ where $r_{i0}=0$ and \vspace{-2mm}
\[
r_i(t)={\bmat{x(t) \\w(t)\\u(t)\\ \dot x(t)}}, \qquad i=1,\cdots,K.\vspace{1mm}
\]
and\vspace{-5mm}
{\small\begin{align*}
&v(t)=\sum_{i=1}^K \bmat{A_i  & B_{1i} & B_{2i}& E_i\\C_{1i}& D_{11i} & D_{12i} & C_{1ei}\\C_{2i} & D_{21i} & D_{22i}&C_{2ei}} \bmat{x(t-\tau_i)\\w(t-\tau_i)\\u(t-\tau_i)\\ \dot x(t-\tau_i)}+\\
&\sum_{i=1}^K \int\limits_{-\tau_i}^0\bmat{A_{di}(s) & B_{1di}(s) &B_{2di}(s)& E_{di}(s)\\C_{1di}(s) & D_{11di}(s) & D_{12di}(s)& C_{1dei}(s)\\C_{2di}(s) &D_{21di}(s) & D_{22di}(s)& C_{2dei}(s)} \bmat{x(t+s)\\w(t+s)\\u(t+s)\\ \dot x(t+s)}ds.
\end{align*}}
\end{lem}
\begin{pf}
Given $u$, $w$, suppose $x$, $y$, and $z$ satisfy the NDS. Defining $r_i$ and $v$, we have
\[
v(t)=\sum_{i=1}^K C_{vi} r_i(t-\tau_i)+\sum_{i=1}^K \int_{-\tau_i}^0C_{vdi}(s) r_i(t+s)ds
\]
as desired. Furthermore,
\[
\bmat{\dot{x}(t)\\ z(t)\\y(t)}=\bmat{A_0 & B_1& B_2\\C_1 &D_{11}&D_{12}\\C_2&D_{21}&D_{22}}\bmat{x(t)\\w(t)\\u(t)}+\bmat{B_v\\D_{1v}\\D_{2v}} v(t)
\]
Now all that remains is to show
\[	
r_i(t)=\bmat{C_{ri}&B_{r1i}&B_{r2i}}\bmat{x(t)\\w(t)\\u(t)}+D_{rvi} v(t)
\]
\begin{align*}
r_i(t)&=\bmat{x(t)\\w(t)\\u(t)\\ \dot x(t) }=\bmat{I & 0& 0\\ 0&I&0\\0&0&I\\0 & 0 & 0}\bmat{x(t)\\w(t)\\u(t)}+\bmat{0\\0\\0\\\dot x(t)}\\
&=\bmat{I & 0& 0\\ 0&I&0\\0&0&I\\A_0 & B_1 & B_2}\bmat{x(t)\\w(t)\\u(t)}+\bmat{0&0&0\\0&0&0\\0&0&0\\I&0&0} v(t)\\
&=\bmat{C_{ri}&B_{r1i}&B_{r2i}}\bmat{x(t)\\w(t)\\u(t)}+D_{rvi} v(t)
\end{align*}
as desired. Finally, $r_{i0}=0$ since $x(t),w(t),u(t)=0$ for $t\le 0$.

\hrulefill

For the converse, given $u$, $w$ suppose $x$, $y$, $z$, $v$ and $r_i$ satisfy the DDF where $r_{i0}=0$. Then \vspace{-2mm}
\begin{align*}
&\bmat{\dot{x}(t)\\ z(t)\\y(t)}=\bmat{A_0 & B_1& B_2\\C_1 &D_{11}&D_{12}\\C_2&D_{21}&D_{22}}\bmat{x(t)\\w(t)\\u(t)}+ v(t)
\end{align*}
and hence
\begin{align*}
&r_i(t)=\bmat{C_{ri}&B_{r1i}&B_{r2i}}\bmat{x(t)\\w(t)\\u(t)}+D_{rvi} v(t) \\
&=\bmat{I & 0& 0\\ 0&I&0\\0&0&I\\A_0 & B_1 & B_2}\bmat{x(t)\\w(t)\\u(t)}+\bmat{0&0&0\\0&0&0\\0&0&0\\I&0&0} v(t) =\bmat{x(t) \\w(t)\\u(t)\\ \dot x(t)}.
\end{align*}
We now have that
\begin{align*}
&v(t)=\sum_{i=1}^K C_{vi} r_i(t-\tau_i)+\sum_{i=1}^K \int_{-\tau_i}^0C_{vdi}(s) r_i(t+s)ds\\
&=\sum_{i=1}^K \bmat{A_i  & B_{1i} & B_{2i}& E_i\\C_{1i}& D_{11i} & D_{12i} & C_{1ei}\\C_{2i} & D_{21i} & D_{22i}&C_{2ei}} r_i(t-\tau_i)+\\
&\sum_{i=1}^K \int\limits_{-\tau_i}^0\bmat{A_{di}(s) & B_{1di}(s) &B_{2di}(s)& E_{di}(s)\\C_{1di}(s) & D_{11di}(s) & D_{12di}(s)& C_{1dei}(s)\\C_{2di}(s) &D_{21di}(s) & D_{22di}(s)& C_{2dei}(s)} r_i(t+s)ds\\
&=\sum_{i=1}^K \bmat{A_i  & B_{1i} & B_{2i}& E_i\\C_{1i}& D_{11i} & D_{12i} & C_{1ei}\\C_{2i} & D_{21i} & D_{22i}&C_{2ei}} \bmat{x(t-\tau_i)\\w(t-\tau_i)\\u(t-\tau_i)\\ \dot x(t-\tau_i)}+\\
&\sum_{i=1}^K \int\limits_{-\tau_i}^0\bmat{A_{di}(s) & B_{1di}(s) &B_{2di}(s)& E_{di}(s)\\C_{1di}(s) & D_{11di}(s) & D_{12di}(s)& C_{1dei}(s)\\C_{2di}(s) &D_{21di}(s) & D_{22di}(s)& C_{2dei}(s)} \bmat{x(t+s)\\w(t+s)\\u(t+s)\\ \dot x(t+s)}ds.
\end{align*}
from which we conclude that
\begin{align*}
&\bmat{\dot{x}(t)\\ z(t)\\y(t)}=\bmat{A_0 & B_1& B_2\\C_1 &D_{11}&D_{12}\\C_2&D_{21}&D_{22}}\bmat{x(t)\\w(t)\\u(t)}+ v(t)\\
&=\bmat{A_0 & B_1 & B_2\\ C_{10} & D_{11} &D_{12}\\ C_{20} & D_{21} &D_{22}}\bmat{x(t)\\w(t)\\u(t)}\notag\\
&+\sum_{i=1}^K \bmat{A_i  & B_{1i} & B_{2i}& E_i\\C_{1i}& D_{11i} & D_{12i} & E_{1i}\\C_{2i} & D_{21i} & D_{22i}&E_{2i}} \bmat{x(t-\tau_i)\\w(t-\tau_i)\\u(t-\tau_i)\\ \dot x(t-\tau_i)}\notag\\
& +\sum_{i=1}^K \hspace{0mm}\int_{-\tau_i}^0\hspace{-1mm}\bmat{A_{di}(s) & \hspace{-1mm}B_{1di}(s) &\hspace{-1mm}B_{2di}(s)& \hspace{-1mm}E_{di}(s)\\C_{1di}(s) & \hspace{-1mm}D_{11di}(s) & \hspace{-1mm}D_{12di}(s)& \hspace{-1mm}E_{1di}(s)\\C_{2di}(s) &\hspace{-1mm}D_{21di}(s) & \hspace{-1mm}D_{22di}(s)& \hspace{-1mm}E_{2di}(s)} \bmat{x(t+s)\\w(t+s)\\u(t+s)\\ \dot x(t+s)}\hspace{-1mm}ds\notag\\[-8mm]\notag
\end{align*}
as desired. Finally, by definition $x_0=0$, which concludes the proof.
\end{pf}

\newpage
\appendix{Proof of \textbf{Lemma 3}}

\begin{lem} Suppose for given $u$, $w$, $r_{i0}$, that $x$, $r_i$, $v$, $y$, and $z$ satisfy the DDF defined by $\{A_{i}, B_i,C_i,D_{ij},\cdots\}$, with $C_{vdi}$ bounded. Then for $u$, $w$, $\phi_{i0}(s)=r_{i0}(\tau_i s)$, we have that $x$, $v$, $y$, and $z$ also satisfy the ODE-PDE defined by $\{A_{i}, B_i,C_i,D_{ij},\cdots\}$ with $
\phi_i(t,s)=r_i(t+\tau_i s)$. Similarly, for given $u$, $w$, $\phi_{i0}$, if $x$, $v$, $y$, $\phi_i$ and $z$ satisfy the ODE-PDE defined by $\{A_{i}, B_i,C_i,D_{ij},\cdots\}$, then $x$, $v$, $y$, and $z$ satisfy the DDF with $r_i(t)=\phi_i(t,0)$ and $r_{i0}(s)=\phi_{i0}(s/\tau_i)$.\vspace{-2mm}
\end{lem}

\begin{pf}
Suppose $x$, $v$, $y$, and $z$ satisfy the DDF, $\phi_i(t,s)=r_i(t+\tau_i s)$ and $\phi_{i0}(s)=r_{i0}(\tau_i s)$. By assumption, $r_{i0} \in W^{1,2}[-\tau_i,0]$ and hence is absolutely continuous with classical derivative $\dot r_{i0}\in L_2[-\tau_i,0]$. Furthermore, $\phi_{i0}\in W^{1,2}[-1,0]$ and is likewise absolutely continuous with classical derivative. Since $r_i(t+\cdot)\in W^{1,2}[-\tau_i,0]$, we have $\phi_i(t,\cdot)=r_i(t+\tau_i \cdot)\in W^{1,2}[-1,0]$ for all $t \ge 0$.

Note that the continuity property may be proven directly and not as an assumption. Specifically, let us examine differentiability of $\phi_i(t,s)=r_i(t+\tau_i s)$. For any $i\in\{1,K\}$, on the interval $t \in [0,\tau_1]$, we have
{\footnotesize \begin{align*}
\dot{r}_i(t)&=C_{ri}\dot{x}(t)+B_{r1i}\dot{w}(t)+B_{r2i}\dot{u}(t)+D_{rvi}\dot{v}(t)\\
&=C_{ri}\left(A_0x(t)+B_1w(t)+B_2u(t)+B_v v(t)\right)\\
&\qquad +B_{r1i}\dot{w}(t)+B_{r2i}\dot{u}(t)+D_{rvi}\dot{v}(t)\\
&=C_{ri}\left(A_0x(t)+B_1w(t)+B_2u(t)+B_v v(t)\right)\\
&\qquad +B_{r1i}\dot{w}(t)+B_{r2i}\dot{u}(t)\\
&\hspace{-6mm}+D_{rvi}\bbbl(\sum_{j=1}^K C_{vj} \dot r_{j0}(t-\tau_j)+\sum_{j=1}^K \int_{t-\tau_j}^{0}\hspace{-2mm}C_{vdj}(s-t) \dot r_{j0}(s)ds\\
&\qquad \qquad +\sum_{j=1}^K \int_{0}^{t}C_{vdj}(s-t) \dot r_{j}(s)ds\bbbr).
\end{align*}}
Now recall $B_{r1i}\dot{w},B_{r2i}\dot{u} \in L_2[0,\infty]$ and $\dot r_{i0}\in L_2[-\tau_i,0]$. Hence, by Gronwall-Bellman, since $C_{vdi}$ are bounded, we conclude that $P_{\tau_1}r_{i} \in W^{1,2}[-\tau_i,\tau_1]$ where $P_{T}$ is the truncation operator at time $T$. Proceeding by steps, we conclude that $P_T r_i\in W^{1,2}[-\tau_i,T]$ for any $T\ge0$. Thus we conclude that $\phi_i(t)\in W^{1,2}[-1,0]$ for all $t \ge 0$.

Now, examining the initial conditions, we find $\phi_{i}(0,s)=r_i(\tau_i s)=r_{i0}(\tau_i s)=\phi_{i0}(s)$ and\vspace{-3mm}
\begin{align*}
&\phi_{i0}(0)=r_{i0}(0)=C_{ri}x_0\\
&+D_{rvi}\left(\sum_{i=1}^K C_{vi} r_{i0}(-\tau_i)+\sum_{i=1}^K \int_{-\tau_i}^0C_{vdi}(s) r_{i0}(s)ds \right)\notag\\
&=C_{ri}x_0\\
&+D_{rvi}\left(\sum_{i=1}^K C_{vi} \phi_{i0}(-1)+\sum_{i=1}^K \int_{-1}^0\tau_i C_{vdi}(\tau_i s) \phi_{i0}(s)ds\right).\notag\\[-8mm] \notag
\end{align*}
Furthermore,\vspace{-3mm}
\begin{align*}
&v(t)=\sum_{i=1}^K C_{vi} r_i(t-\tau_i)+\sum_{i=1}^K \int_{-\tau_i}^0C_{vdi}( s) r_i(t-\tau_i s)ds\\
&=\sum_{i=1}^K C_{vi} \phi_i(t,-1)+\sum_{i=1}^K \int_{-1}^0\tau_iC_{vdi}(\tau_i s) \phi_i(t, s)ds
\end{align*}
and\vspace{-3mm}
\begin{align*}
\dot \phi_i(t,s)&=\dot r_i(t+\tau_i s)=\frac{1}{\tau_i}\partial_sr_i(t+\tau_i s)=\frac{1}{\tau_i}\phi_{i,s}(t,s).
\end{align*}
Finally, $x$, $z$, $y$ satisfy the ODE-PDE by inspection and\vspace{-3mm}
\[
\phi_i(t,0)=r_i(t)=C_{ri}x(t)+B_{r1i}w(t)+B_{r2i}u(t)+D_{rvi}v(t).\vspace{-3mm}
\]
\hrulefill\vspace{-3mm}

Conversely, suppose $x$, $v$, $y$, $\phi_i$ and $z$ satisfy the ODE-PDE and $r_{i0}(s)=\phi_{i0}(s/\tau_i)$. Since $\phi_i(t)\in W^{1,2}[-1,0]$, $r_i(t+s)=\phi_i(t,s/\tau_i)$ and hence $r_i(t+\cdot)\in W^{1,2}[-\tau_i,0]$. Likewise, $r_{i0}\in W^{1,2}[-\tau_i,0]$, and\vspace{-3mm}
\begin{align*}
&r_{i0}(0)=\phi_{i0}(0)
=C_{ri}x_0\\
&+D_{rvi}\left(\sum_{i=1}^K C_{vi} \phi_{i0}(-1)+\sum_{i=1}^K \int_{-1}^0\tau_i C_{vdi}(\tau_i s) \phi_{i0}(s)ds\right)\notag\\
&=C_{ri}x_0\\
&+D_{rvi}\left(\sum_{i=1}^K C_{vi} r_{i0}(-\tau_i)+\sum_{i=1}^K \int_{-\tau_i}^0C_{vdi}(s) r_{i0}(s)ds \right).
\end{align*}
Now if $r_i(t)=\phi_i(t,0)$, then\vspace{-3mm}
\[
r_i(t)=\phi_i(t,0)=C_{ri}x(t)+B_{r1i}w(t)+B_{r2i}u(t)+D_{rvi}v(t).\vspace{-3mm}
\]
Next, if $\dot \phi_i(t,s)=\frac{1}{\tau_i}\partial_s \phi_{i}(t,s)$, then $\phi_i(t,s)=\phi_i(t+s\tau_i,0)=r_i(t+s\tau_i)$ and hence\vspace{-3mm}
\begin{align*}
&v(t)\\
&=\sum_{i=1}^K C_{vi} \phi_i(t,-1)+\sum_{i=1}^K \int_{-1}^0\tau_iC_{vdi}(\tau_i s) \phi_i(t, s)ds\\
&=\sum_{i=1}^K C_{vi} r_i(t-\tau_i)+\sum_{i=1}^K \int_{-\tau_i}^0C_{vdi}( s) r_i(t-\tau_i s)ds.
\end{align*}
Finally, $x$, $z$, $y$ satisfy the DDF by inspection and for $s\in[-\tau_i,0]$, $r_i(s)=\phi_i(0,s/\tau_i)=\phi_{i0}(s/\tau_i)=r_{i0}(s)$.
%
%\begin{align*}
%\dot r_i(t+\tau_i s) = \dot \phi_i(t,s)&=\frac{1}{\tau_i}\partial_s(t+\tau_i s)=\frac{1}{\tau_i}\phi_{i,s}(t,s)\notag\\
%\end{align*}
\end{pf}

\newpage

\appendix{Proof of \textbf{Lemma 4}}

\begin{lem}[Lemma 4]
  Given $u$, $w$, and $x_0$, $\phi_{i0}$ satisfying Condition~\eqref{eqn:sewing_PDE},  Suppose $x$, $\phi_i$, $v$, $y$, and $z$ satisfy the ODE-PDE defined by $\{A_{i}, B_i,C_i,D_{ij},\cdots\}$. Then $y$ and $z$ also satisfy the PIE defined by $\{\mcl A, \mcl B_i, \mcl C_i, \mcl D_{ij},\mcl B_{T_i}\}$ with $\mcl T, \mcl A, \mcl B_{i}, \mcl C_i, \mcl D_{ij}$ as defined in Eqn.~\eqref{eqn:PIE_ops} and\vspace{-2mm}
\[
\mbf x(t):={\scriptsize\bmat{x(t)\\ \partial_s \phi_{1}(t,\cdot)\\\vdots\\ \partial_s \phi_{K}(t,\cdot)}}\quad \mbf x_0:={\scriptsize\bmat{x_0\\ \partial_s \phi_{10}\\ \vdots \\ \partial_s \phi_{K0}}}.\vspace{-2mm}
\]
Furthermore, for given $u$, $w$, $\mbf x_{0} \in \R^n \times L_2[-1,0]^p$, if $y$, $z$ and $\mbf x$ satisfy the PIE defined by $\{\mcl A, \mcl B_i, \mcl C_i, \mcl D_{ij},\mcl B_{T_i}\}$, then %for
%\[
%\bmat{x_0\\ \phi_{10}\\\vdots\\ \phi_{K0}}=\mcl T \mbf x_0
%\]
 $x$, $\phi_i$, $v$, $y$, and $z$ satisfy the ODE-PDE defined by $\{A_{i}, B_i,C_i,D_{ij},\cdots\}$ where\vspace{-2mm}
\[
{\scriptsize\bmat{x(t)\\ \phi_{1}(t,\cdot)\\\vdots\\ \phi_{K}(t,\cdot)}}=\mcl T \mbf x(t)+\mcl B_{T1}w(t)+\mcl B_{T2}u(t),\; {\scriptsize\bmat{x_0\\ \phi_{10}\\\vdots\\ \phi_{K0}}}=\mcl T \mbf x_0.\vspace{-2mm}
\]
\end{lem}
\begin{pf}Suppose $x$, $\phi_i$, $v$, $y$, and $z$ satisfy the ODE-PDE and\vspace{-3mm}
\[
\mbf x(t):=\bmat{x(t)\\ \partial_s \phi_{1}(t,\cdot)\\\vdots\\ \partial_s \phi_{K}(t,\cdot)}\quad \mbf x_0:=\bmat{x_0\\ \partial_s \phi_{10}\\ \vdots \\ \partial_s \phi_{K0}}.\vspace{-2mm}
\]
Then, $\mbf x_0\in \R^n \times L_2[-1,0]^p$ and
\begin{equation}
\phi_i(t,s)=\phi_i(t,0)-\int_s^0  \phi_{i,s}(t,\eta)d \eta.\label{eqn:FS}
\end{equation}
%To obtain the PIE representation, we eliminate $\phi_i$, $\phi_i(0)$, and $\phi_i(-1)$ from Eqns.~\eqref{eqn:ODEPDE}.
%
%The main challenge is to solve for $v$ in terms of $x$, $w$, $u$ and $\phi_{i,s}$. Specifically, we seek to show that
%\begin{align*}
%v(t)=&C_{vx}x(t)+D_{vw}w(t)+D_{vu}u(t)\\
%&\qquad  +\int_{-1}^0  \sum_{i=1}^K C_{Ii}(s)  \phi_{i,s}(t,s)d s.
%\end{align*}
Recall from Eqns.~\eqref{eqn:ODEPDE} that\vspace{-3mm}
\[
v(t)=\sum_{i=1}^K C_{vi} \phi_i(t,-1)+\sum_{i=1}^K \int_{-1}^0\tau_iC_{vdi}(\tau_i s) \phi_i(t, s)ds\vspace{-3mm}
\]
and\vspace{-3mm}
\[
\phi_i(t,0)=C_{ri}x(t)+B_{r1i}w(t)+B_{r2i}u(t)+D_{rvi}v(t).\vspace{-3mm}
\]
From Eqn.~\eqref{eqn:FS},\vspace{-3mm}
\[
\phi_i(t,-1)=\phi_i(t,0) -\int_{-1}^0 \mbf \phi_{i,s}(t,\eta)d \eta\vspace{-3mm}
\]
and hence\vspace{-3mm}
\begin{align*}
v(t)&=\sum_{i=1}^K C_{vi} \left(\phi_i(t,0)-\int_{-1}^0 \mbf \phi_{i,s}(t,\eta)d \eta\right)\\
&\hspace{2mm} +\sum_{i=1}^K \int_{-1}^0\tau_iC_{vdi}(\tau_i s) \left(\phi_i(t,0)-\int_{s}^0 \mbf \phi_{i,s}(t,\eta)d \eta\right)ds\\
&=\left(\sum_{i=1}^K \hat C_{vi}C_{ri}\right)x(t)+\left(\sum_{i=1}^K \hat C_{vi}B_{r1i}\right)w(t)\\
&\hspace{2mm}+\left(\sum_{i=1}^K \hat C_{vi}B_{r2i}\right)u(t)+\left(\sum_{i=1}^K \hat C_{vi}D_{rvi}\right)v(t)\\
&\hspace{2mm}-\left( \int_{-1}^0 \sum_{i=1}^K C_{vi} \mbf \phi_{i,s}(t,\eta)d \eta\right)\\
&\hspace{2mm} -\left(\sum_{i=1}^K \tau_i \int_{-1}^0\int_{s}^0 C_{vdi}(\tau_i s) \mbf \phi_{i,s}(t,\eta)d \eta ds\right).
\end{align*}
Eliminating $v$ from the RHS, we obtain\vspace{-3mm}
\begin{align*}
v(t)=&C_{vx}x(t)+D_{vw}w(t)+D_{vu}u(t)\\
&-D_I\left( \int_{-1}^0 \sum_{i=1}^K C_{vi} \mbf \phi_{i,s}(t,\eta)d \eta\right)\\
& -D_I\left(\sum_{i=1}^K \tau_i \int_{-1}^0\int_{s}^0 C_{vdi}(\tau_i s) \mbf \phi_{i,s}(t,\eta)d \eta ds\right).
\end{align*}
Using the identity\vspace{-3mm}
\[
 \int_{-1}^0\int_{s}^0f(s,\eta)d\eta ds=\int_{-1}^0\int_{-1}^s f(\eta,s)  d\eta ds,\vspace{-3mm}
\]
we obtain\vspace{-3mm}
\begin{align*}
&v(t)=C_{vx}x(t)+D_{vw}w(t)+D_{vu}u(t)\\
&\qquad-D_I\left( \int_{-1}^0 \sum_{i=1}^K C_{vi} \mbf \phi_{i,s}(t,s)d s\right) \\
&\qquad-D_I\left(\sum_{i=1}^K  \int_{-1}^0\left( \int_{-1}^s \tau_iC_{vdi}(\tau_i \eta) d \eta\right)\mbf \phi_{i,s}(t,s) ds\right)\\
&=C_{vx}x(t)+D_{vw}w(t)+D_{vu}u(t)\\
&-D_I\left( \int_{-1}^0 \sum_{i=1}^K \left(C_{vi} +\tau_i\int_{-1}^sC_{vdi}(\tau_i \eta) d \eta \right) \mbf \phi_{i,s}(t,s)d s\right)\\
&=C_{vx}x(t)+D_{vw}w(t)+D_{vu}u(t)\\
&\qquad + \int_{-1}^0  \sum_{i=1}^K C_{Ii}(s) \mbf \phi_{i,s}(t,s)d s.
\end{align*}
The rest of the sufficiency proof is straightforward. Plugging this expression for $v(t)$ into Eqns.~\eqref{eqn:ODEPDE}, we obtain\vspace{-3mm}
\begin{align*}
&z(t)=\mbf C_{10}x(t)+\int\limits_{-1}^0 \mbf C_{11}(s)\mbf \Phi(t,\eta)d \eta+\mbf D_{11}w(t)+\mbf D_{12}u(t)\\
&=\mcl C_1\mbf x(t)+\mcl D_{11}w(t)+\mcl D_{12}u(t)\\
&y(t)=\mbf C_{20}x(t)+\int\limits_{-1}^0 \mbf C_{21}(s)\mbf \Phi(t,\eta)d \eta+\mbf D_{21}w(t)+\mbf D_{22}u(t)\\
&=\mcl C_2\mbf x(t)+\mcl D_{21}w(t)+\mcl D_{22}u(t)
\end{align*}
where\vspace{-3mm}
\[
\mbf x(t):=\bmat{x(t)\\ \mbf \Phi(t,\cdot)},\quad \mbf \Phi(t):=\bmat{ \phi_{1,s}(t,\cdot)\\\vdots\\ \phi_{K,s}(t,\cdot)}.\vspace{-3mm}
\]
Likewise,\vspace{-3mm}
\begin{align*}
&\bmat{\dot x(t)\\ \dot \phi_1(t,s)\\\vdots \\\dot \phi_K(t,s)}=\bmat{\mbf A_0 x(t)+\int\limits_{-1}^0 \mbf A(\eta) \mbf \Phi(t,\eta)d \eta\\I_\tau \mbf \Phi(t,s)}\\
&\qquad +\bmat{\mbf B_1 \\0}w(t)+\bmat{\mbf B_2\\0}u(t)\\
&=\mcl A\mbf x(t)+\mcl B_1w(t)+\mcl B_2u(t).
\end{align*}
Finally, we observe that\vspace{-3mm}
\begin{align*}
&\phi_i(t,s)=C_{ri}x(t)+B_{r1i}w(t)+B_{r2i}u(t)+D_{rvi}v(t)\\
& \qquad\qquad\qquad-\int_{s}^0 \mbf \phi_{i,s}(t,\eta)d \eta\\
&=(C_{ri}+D_{rvi}C_{vx}){x}(t)+(B_{r1i}+D_{rvi}D_{vw}){w}(t)\\
&\qquad\qquad\qquad+(B_{r2i}+D_{rvi}D_{vu}){u}(t)\\
&-\int_{s}^0 {\mbf\phi}_{i,s}(t,\eta)d \eta+\left( \int_{-1}^0  \sum_{j=1}^K D_{rvi} C_{Ij}(s) {\mbf \phi}_{j,s}(t,s)d s\right).
\end{align*}
Hence\vspace{-3mm}
\begin{align*}
&\bmat{\phi_{1}(t,\cdot)\\\vdots\\ \phi_{K}(t,\cdot)}=\mbf T_0x(t)+\mbf T_1{w}(t)+\mbf T_2{u}(t)\\
&+\int_{-1}^s \mbf T_a(\eta)\mbf\Phi(t,\eta)d \eta+\int_{s}^0 \mbf T_b(\eta)\mbf\Phi(t,\eta)d \eta
\end{align*}
and\vspace{-3mm}
\[
\bmat{x(t)\\ \phi_{1}(t,\cdot)\\\vdots\\ \phi_{K}(t,\cdot)}=\mcl T \mbf x(t)+\mcl B_{T1}w(t)+\mcl B_{T2}u(t).\vspace{-3mm}
\] Finally, we differentiate and combine these expressions to obtain\vspace{-3mm}
%\begin{align*}
%\dot{\mbf x}_p(t)&=\mcl A\mbf x(t)+\mcl B_1w(t)+\mcl B_2u(t)\\
%z(t)&=\mcl C_1\mbf x(t)+\mcl D_{11}w(t)+\mcl D_{12}u(t),\\
%y(t)&=\mcl C_2\mbf x(t)+\mcl D_{21}w(t)+\mcl D_{22}u(t).
%\end{align*}
\begin{align*}
\mcl T \dot{\mbf x}(t)+\mcl B_{T_1}\dot w(t)+\mcl B_{T_2}\dot u(t)&=\mcl A\mbf x(t)+\mcl B_1w(t)+\mcl B_2u(t)\notag\\
\hspace{-1cm}z(t)=\mcl C_1\mbf x(t)&+\mcl D_{11}w(t)+\mcl D_{12}u(t),\notag\\
\hspace{-1cm}y(t)=\mcl C_2\mbf x(t)&+\mcl D_{21}w(t)+\mcl D_{22}u(t).
\end{align*}
We conclude that $\mbf x$, $y$ and $z$ satisfy the PIE.\vspace{-3mm}

\hrulefill

Conversely, suppose $\mbf x$, $y$ and $z$ satisfy the PIE. Partition $\mbf x$ as\vspace{-3mm}
\[
\mbf x(t,s)=\bmat{x_1(t)\\\Phi(t,s)}\vspace{-3mm}
\]
 and define\vspace{-3mm}
\[
{\scriptsize\bmat{x(t)\\ \phi_{1}(t,\cdot)\\\vdots\\ \phi_{K}(t,\cdot)}}
=\mcl T \mbf x(t)+\mcl B_{T1}w(t)+\mcl B_{T2}u(t),\; {\scriptsize\bmat{x_0\\ \phi_{10}\\\vdots\\ \phi_{K0}}}=\mcl T \mbf x_0,\vspace{-2mm}
\]
and\vspace{-3mm}
\[
v(t)=C_{vx}x(t)+D_{vw}w(t)+D_{vu}u(t)
+\int\limits_{-1}^0\sum_{j=1}^K C_{Ij}\Phi_{j}(t,\theta)d\theta.\vspace{-3mm}
\]
Examine\vspace{-3mm}
\[
\mcl T \mbf x(t)=\bmat{x(t)\\T_0x(t)+\int_{-1}^0T_a(\theta)\Phi(t,\theta)d\theta - \int_{s}^0\Phi(t,\theta)d\theta}.\vspace{-3mm}
\]
This implies that\vspace{-3mm}
\[
\bmat{\phi_{1,s}(t,\cdot)\\\vdots\\ \phi_{K,s}(t,\cdot)}=\Phi(t,s)\vspace{-3mm}
\]
and\vspace{-3mm}
\begin{align*}
&\bmat{x(t)\\\phi_{1}(t,0)\\\vdots\\ \phi_{K}(t,0)}\\
&=\bmat{x_1(t)\\T_0x(t)+\int\limits_{-1}^0T_a(\theta)\Phi(t,\theta)d\theta}+\mcl B_{T1}w(t)+\mcl B_{T2}u(t)\\
\end{align*}
and hence $x(t)=x_1(t)$ and\vspace{-3mm}
{\scriptsize \begin{align*}
&\phi_i(t,0)=C_{ri}x(t) +B_{r1i}w(t)+B_{r2i}u(t)\\
&+
D_{rvi}\left(C_{vx}x(t)+D_{vw}w(t)+D_{vu}u(t)
+\int\limits_{-1}^0\sum_{j=1}^K C_{Ij}\phi_{j,s}(t,\theta)d\theta\right)\\
&=C_{ri}x(t) +B_{r1i}w(t)+B_{r2i}u(t)+D_{rvi}v(t)
\end{align*}}
as desired. Next, we have\vspace{-3mm}
\begin{align*}
&{\scriptsize\bmat{\dot x(t)\\ \dot \phi_{1}(t,\cdot)\\\vdots\\ \dot \phi_{K}(t,\cdot)}}
=\mcl T \dot{\mbf x}(t)+\mcl B_{T1}\dot w(t)+\mcl B_{T2}\dot u(t)\\
&=\mcl A\mbf x(t)+\mcl B_1w(t)+\mcl B_2u(t)\\
&=\bmat{A_0x(t) +B_1w(t)+B_2u(t)\\\frac{1}{\tau_1}\phi_{1,s}(t,s)\\\vdots \\ \frac{1}{\tau_K}\phi_{K,s}(t,s)}\\
&+\scriptsize{\bmat{B_v \left(C_{vx}x(t)+D_{vw}w(t)+D_{vu}u(t)+\sum\nolimits_{j=1}^K C_{Ij}\phi_{j,s}(t,\theta)d\theta \right)\\0}}\\
&=\bmat{A_0x(t) +B_1w(t)+B_2u(t)+B_v v(t)\\\frac{1}{\tau_1}\phi_{1,s}(t,s)\\\vdots \\ \frac{1}{\tau_K}\phi_{K,s}(t,s)}
\end{align*}
as desired. Returning to $v$, we have\vspace{-3mm}
%\[
%v(t)=C_{vx}x(t)+D_{vw}w(t)+D_{vu}u(t)
%+\int\limits_{-1}^0\sum_{j=1}^K C_{Ij}\Phi_{j}(t,\theta)d\theta
%\]
\begin{align*}
&v(t)\\
&=C_{vx}x(t)+D_{vw}w(t)+D_{vu}u(t)
+\int\limits_{-1}^0\sum_{j=1}^K C_{Ij}\Phi_{j}(t,\theta)d\theta\\
&=C_{vx}x(t)+D_{vw}w(t)+D_{vu}u(t)\\
&-D_I\left( \int_{-1}^0 \sum_{i=1}^K \left(C_{vi} +\tau_i\int_{-1}^sC_{vdi}(\tau_i \eta) d \eta \right) \mbf \phi_{i,s}(t,s)d s\right)\\
%&=C_{vx}x(t)+D_{vw}w(t)+D_{vu}u(t)\\
%&\qquad-D_I\left( \int_{-1}^0 \sum_{i=1}^K C_{vi} \mbf \phi_{i,s}(t,s)d s\right) \\
%&\qquad-D_I\left(\sum_{i=1}^K  \int_{-1}^0\left( \int_{-1}^s \tau_iC_{vdi}(\tau_i \eta) d \eta\right)\mbf \phi_{i,s}(t,s) ds\right)\\
%&=C_{vx}x(t)+D_{vw}w(t)+D_{vu}u(t)\\
%&-D_I\left( \int_{-1}^0 \sum_{i=1}^K C_{vi} \mbf \phi_{i,s}(t,\eta)d \eta\right)\\
%& -D_I\left(\sum_{i=1}^K \tau_i \int_{-1}^0\int_{s}^0 C_{vdi}(\tau_i s) \mbf \phi_{i,s}(t,\eta)d \eta ds\right)\\
&=D_I\bbbbl(\left(\sum_{i=1}^K \hat C_{vi}C_{ri}\right)x(t)+\left(\sum_{i=1}^K \hat C_{vi}B_{r1i}\right)w(t)\\
&+\left(\sum_{i=1}^K \hat C_{vi}B_{r2i}\right)u(t)-\left( \int_{-1}^0 \sum_{i=1}^K C_{vi} \mbf \phi_{i,s}(t,\eta)d \eta\right)\\
& -\left(\sum_{i=1}^K \tau_i \int_{-1}^0\int_{s}^0 C_{vdi}(\tau_i s) \mbf \phi_{i,s}(t,\eta)d \eta ds\right)\bbbbr)\\
&=D_I\bbbbl(\sum_{i=1}^K \hat C_{vi}\left(C_{ri}x(t)+B_{r1i}w(t)+B_{r2i}u(t)\right)\\
&\qquad -\sum_{i=1}^K C_{vi} \left(\int_{-1}^0   \mbf \phi_{i,s}(t,\eta)d \eta\right)\\
&\qquad  -\sum_{i=1}^K  \int_{-1}^0\tau_i \left( \int_{s}^0 C_{vdi}(\tau_i s) \mbf \phi_{i,s}(t,\eta)d \eta \right)ds\bbbbr).
\end{align*}
Hence\vspace{-3mm}
\begin{align*}
&\left(I-\left(\sum_{i=1}^K \hat C_{vi}D_{rvi}\right)\right)v(t)=\\
&=\sum_{i=1}^K \hat C_{vi}\left(C_{ri}x(t)+B_{r1i}w(t)+B_{r2i}u(t)\right)\\
&-\sum_{i=1}^K C_{vi} \left(\int_{-1}^0   \mbf \phi_{i,s}(t,\eta)d \eta\right)\\
& -\sum_{i=1}^K  \int_{-1}^0\tau_i \left( \int_{s}^0 C_{vdi}(\tau_i s) \mbf \phi_{i,s}(t,\eta)d \eta \right)ds.
\end{align*}
Hence\vspace{-3mm}
\begin{align*}
&v(t)=\\
&=\sum_{i=1}^K \hat C_{vi}\left(C_{ri}x(t)+B_{r1i}w(t)+B_{r2i}u(t)+D_{rvi}v(t)\right)\\
&\qquad -\sum_{i=1}^K C_{vi} \left(\int_{-1}^0   \mbf \phi_{i,s}(t,\eta)d \eta\right)\\
&\qquad -\sum_{i=1}^K  \int_{-1}^0\tau_i \left( \int_{s}^0 C_{vdi}(\tau_i s) \mbf \phi_{i,s}(t,\eta)d \eta \right)ds\\
&=\sum_{i=1}^K \hat C_{vi}\phi_i(t,0)-\sum_{i=1}^K C_{vi} \left(\int_{-1}^0   \mbf \phi_{i,s}(t,\eta)d \eta\right)\\
&\qquad -\sum_{i=1}^K  \int_{-1}^0\tau_i \left( \int_{s}^0 C_{vdi}(\tau_i s) \mbf \phi_{i,s}(t,\eta)d \eta\right)ds\\
&=\sum_{i=1}^K \left(C_{vi} +\int_{-1}^0 \tau_iC_{vdi}(\tau_i s)ds\right)\phi_i(t,0)\\
&\qquad-\sum_{i=1}^K C_{vi} \left(\int_{-1}^0   \mbf \phi_{i,s}(t,\eta)d \eta\right)\\
&\qquad -\sum_{i=1}^K  \int_{-1}^0\tau_i \left( \int_{s}^0 C_{vdi}(\tau_i s) \mbf \phi_{i,s}(t,\eta)d \eta\right)ds\\
&=\sum_{i=1}^K C_{vi} \left(\phi_i(t,0) -\int_{-1}^0   \mbf \phi_{i,s}(t,\eta)d \eta\right)\\
& +\sum_{i=1}^K \left(\tau_i \int_{-1}^0 C_{vdi}(\tau_i s) \left(  \phi_i(t,0) -\int_{s}^0\mbf \phi_{i,s}(t,\eta)d \eta\right) ds\right)\\
&=\sum_{i=1}^K C_{vi}\phi_i(t,-1) +\sum_{i=1}^K \left(\tau_i \int_{-1}^0 C_{vdi}(\tau_i s) \phi_i(t,s) ds\right)\end{align*}
as desired. Next, we have
{\scriptsize\begin{align*}
&z(t)=\mcl C_1\mbf x(t)+\mcl D_{11}w(t)+\mcl D_{12}u(t)\\
&=\mbf C_{10}x(t)+\int\limits_{-1}^0 \mbf C_{11}(s)\mbf \phi(t,\eta)d \eta+\mbf D_{11}w(t)+\mbf D_{12}u(t)\\
&=C_{1}x(t)+D_{11}w(t)+D_{12}u(t)\\
&+D_{1v}\left(C_{vx}x(t)+D_{vw}w(t)+D_{vu}u(t)+\int\limits_{-1}^0 \sum_{i=1}^K C_{1i}(s) \phi_{i,s}(t,\eta)d \eta\right)\\
&=C_{1}x(t)+D_{11}w(t)+D_{12}u(t)+D_{1v}v(t)
\end{align*}}
as desired. Furthermore,\vspace{-3mm}
\begin{align*}
&y(t)=\mcl C_2\mbf x(t)+\mcl D_{21}w(t)+\mcl D_{22}u(t)
\\
&=\mbf C_{20}x(t)+\int\limits_{-1}^0 \mbf C_{21}(s)\mbf \Phi(t,\eta)d \eta+\mbf D_{21}w(t)+\mbf D_{22}u(t)\\
&=C_{2}x(t)+D_{21}w(t)+D_{22}u(t)\\
&+D_{2v}\bbbl(C_{vx}x(t)+D_{vw}w(t)+D_{vu}u(t)\\
&\qquad\qquad \qquad +\int\limits_{-1}^0 \sum_{i=1}^K C_{1i}(s) \phi_{i,s}(t,\eta)d \eta\bbbr)\\
&=C_{2}x(t)+D_{21}w(t)+D_{22}u(t)+D_{2v}v(t)
\end{align*}
as desired.
Finally, since\vspace{-3mm}
\[
\bmat{x_0\\ \phi_{10}\\\vdots\\ \phi_{K0}}=\mcl T \mbf x_0.\vspace{-2mm}
\]
as before, we have\vspace{-3mm}
\[
\quad {\scriptsize\bmat{x_0\\ \partial_s \phi_{10}\\ \vdots \\ \partial_s \phi_{K0}}}=\mbf x_0=\bmat{x_0\\\mbf \Phi_0}.\vspace{-2mm}
\]
and\vspace{-3mm}
\begin{align*}
&\bmat{x_0\\\phi_{10}(t,0)\\\vdots\\ \phi_{K0}(t,0)}\\
&=\bmat{x(t)\\T_0x(t)+\int\limits_{-1}^0T_a(\theta)\Phi(t,\theta)d\theta}+\mcl B_{T1}w(t)+\mcl B_{T2}u(t)\\
\end{align*}
and hence\vspace{-3mm}
\begin{align*}
&\phi_{i0}(0)=C_{ri}x_0 \\
&+
D_{rvi}\left(C_{vx}x_0+\int\limits_{-1}^0\sum_{j=1}^K C_{Ij}\phi_{j0,s}(\theta)d\theta\right)\\
&=C_{ri}x_0 +D_{rvi}v(0)\\
&=C_{ri}x_0\\
&+D_{rvi}\left(\sum_{i=1}^K C_{vi} \phi_{i0}(-1)+\sum_{i=1}^K \int_{-1}^0\tau_iC_{vdi}(\tau_i s) \phi_{i0}(s)ds\right)
\end{align*}
as desired.

%\begin{align*}
%&\phi_i(t,0)=C_{ri}x(t) +B_{r1i}w(t)+B_{r2i}u(t)\\
%&+
%D_{rvi}\left(C_{vx}x(t)+D_{vw}w(t)+D_{vu}u(t)
%+\int\limits_{-1}^0\sum_{j=1}^K C_{Ij}\phi_{j,s}(t,\theta)d\theta\right)\\
%&=C_{ri}x(t) +B_{r1i}w(t)+B_{r2i}u(t)+D_{rvi}v(t)
%\end{align*}
%\begin{align*}
%&\bmat{x_0\\\phi_{1}(t,0)\\\vdots\\ \phi_{K}(t,0)}\\
%&=\bmat{x_0\\T_0x_0+\int\limits_{-1}^0T_a(\theta)\Phi_0(t,\theta)d\theta}\\
%\end{align*}
%and hence
%\begin{align*}
%&\phi_{i0}(0)=C_{ri}x_0 \\
%&+
%D_{rvi}\left(C_{vx}x_0+\int\limits_{-1}^0\sum_{j=1}^K C_{Ij}\phi_{j0,s}(\theta)d\theta\right)\\
%&=C_{ri}x_0 +D_{rvi}v(0)\\
%&=\sum_{i=1}^K C_{vi}\phi_{i0}(-1) +\sum_{i=1}^K \left(\tau_i \int_{-1}^0 C_{vdi}(\tau_i s) \phi_{i0}(s) ds\right)
%\end{align*}
%as desired.
\end{pf}

\newpage

\appendix{Proof of \textbf{Lemma 5}}

\begin{lem}[Lemma 5]
For given $r_{i0}$, $x_0$, suppose $r_i$, $v$, $y$, $x$, and $z$ satisfy the DDF defined by \vspace{-2mm}
\[
\{A_0,B_1,B_v,\tilde C_{1},\tilde D_{12},D_{1v},D_{2v},C_{ri},B_{r1i},D_{rvi},C_{vi}\}\vspace{-2mm}
\] given by Eqns.~\eqref{eqn:UAV_SOF_defs}. Then $x$, $z$ and $y$ also satisfy Eqns.~\eqref{eqn:UAV_SOF}.
\end{lem}
\begin{pf}
Suppose $r_i$, $v$, $y$, $x$, and $z$ satisfy the DDF in Eqns.~\eqref{eqn:DDF}. Then\vspace{-3mm}
\begin{align*}
y(t)&=C_{2} x(t)+D_{21}w(t)+D_{2v}v(t)\\
&=C_{2} x(t)+D_{21}w(t)+\sum\limits_{i=1}^N D_{22i}v_i(t)
%&=C_{2} x(t)+D_{21}w(t)+\sum_i D_{22i}r(t-\tau_i)
\end{align*}
and hence\vspace{-3mm}
\begin{align*}
r_i(t)&=C_{ri} x(t)+B_{r1i}w(t)+D_{rvi}v(t)\\
&=FC_{2} x(t)+FD_{21}w(t)+\sum\limits_{i=1}^N FD_{22i}v_i(t)\\
&=F\left(C_{2} x(t)+D_{21}w(t)+\sum_{i=1}^N D_{22i}v_i(t)\right)=Fy(t).
\end{align*}
Next,\vspace{-1mm}
\[
v_i(t)=r_i(t-\tau_i)=Fy(t-\tau_i)\vspace{-1mm}
\]
and we conclude\vspace{-3mm}
\[
y(t)=C_{2} x(t)+D_{21}w(t)+\sum_{i=1}^N D_{22i}Fy(t-\tau_i).\vspace{-3mm}
\]
Similarly\vspace{-3mm}
\begin{align*}
\dot x(t)&=A_0x(t)+B_{1}w(t)+B_{v} v(t)\\
&=A_0x(t)+B_{1}w(t)+\sum\limits_{i=1}^NB_{2i}v_i(t)\\
&=A_0x(t)+B_{1}w(t)+\sum\limits_{i=1}^NB_{2i}Fy(t-\tau_i)
\end{align*}
and finally,\vspace{-3mm}
\begin{align*}
&z(t)=\tilde C_1x(t)+\tilde D_{12}w(t)+D_{1v}v(t)\\
%&=\left(C_1+D_{12}F C_2\right)x(t)+D_{12}FD_{21}w(t)+\sum_i D_{12}F D_{22i}v_i(t)\\
&=C_1x(t)+D_{12}F\left(C_{2} x(t)+D_{21}w(t)+\sum\limits_{i=1}^N D_{22i}v_i(t)\right)\\
&=C_1x(t)+D_{12}Fy(t)
\end{align*}
as desired.
\end{pf}

\end{document}